\documentclass[a4paper]{article}
\usepackage[utf8]{inputenc}
\usepackage[a4paper,left=3cm,right=3cm,top=2cm,bottom=3cm,bindingoffset=5mm]{geometry}
\usepackage[english]{babel}
\usepackage[T1]{fontenc}
\usepackage{graphicx, subfig}
\graphicspath{{img/}}
\usepackage{fancyhdr}
\usepackage{lmodern}
\usepackage{color}

\usepackage{amsfonts}
\usepackage{amsmath}
\usepackage{hyperref}

\title{High symmetric $n$-polygons}
\author{Rolf Haag \\
\href{rhaag.98@gmail.com}{rhaag.98@gmail.com}}
\date{}

\begin{document}
\maketitle

\begin{abstract}
The present article includes the enumeration of $n$-polygons with a certain symmetry property: For an even number $n$ of vertices, we count the $n$-polygons with $\frac{n}{2}$ symmetry axes. In addition, if $n$ is a power of 2, we show the relation to the perfect numbers.\\
\begin{center}
\textbf{Acknowledgments}
\end{center}
I thank Prof. Dr. Hanspeter Kraft from the university of Basel, who supported me in my research by giving me valuable information on form and content and by submitting a template for the basic definitions.
\begin{center}
\textbf{Keywords}
\end{center}
Hamiltonian cycles $\cdot$ Perfect numbers $\cdot$ Polygons $\cdot$ Symmetry-classes
\end{abstract}

\tableofcontents

\section{Introduction}
\label{sec:introduction}
\subsection{Definition of a $n$-polygon}
\label{subsec:definition_of_a_n_polygon}
$n$ vertices are regularly distributed in a circle. We consider the Hamiltonian cycles through the $n$ vertices. \cite{Herman2019} In this paper such Hamiltonian cycles are called $n$-polygons. The usual polygons are the special case where all edges have minimal length.\\

Let $ n $ be a natural number $ n\geq3 $ and $ S ^ 1 \subset \mathbb{R} ^ 2 = \mathbb{C} $ the unit circle in the Euclidean plane. The finite subset
\begin{center}
$V_n:=\{v_k:=e^{2 \pi i k/n} \mid k = 0,1,\ldots,n-1 \} \subset S^1$
\end{center}
represents the vertices of an $n$-polygon. \\

To describe the $n$ polygons we use the $n$-cycles $\sigma = (\sigma_1, \sigma_2, \cdots, \sigma_n)$ consisting of the $n$ numbers $\lbrace 0,1, \cdots , n-1 \rbrace$ in any order. The associated $n$-polygon $P(\sigma)$ is given by the path $\overline{v_{\sigma_1}v_{\sigma_2} \cdots v_{\sigma_n}v_{\sigma_1}}$ or more precisely by combining the links $\overline{v_{\sigma_i}v_{\sigma_{i+1}}}$, $i= 1,2, \cdots , n$, where $\sigma_{n+1}= \sigma_1$. \\

Each of the n edges is assigned its ``length'' $e_i$. $e_i = 1$ means that the $i$-th edge runs counterclockwise from the vertex $V_i$ to the following vertex $V_{i + 1}$. $e_i = 2$ means that the $i$-th edge runs counterclockwise from the vertex $V_i$  to the vertex $V_{i + 2}$ and so on. $e_i = n$ is not possible, since this would mean the connection of the vertex $V_i$ to itself. Therefore, only the numbers between $1$ and $n-1$ are allowed to describe the ``length'' of the edges. The ``length'' of an edge $e_i$ is referred to briefly as a side of the n-polygon. Therefore, an n-polygon can also get described by the $n$-cycle of its sides:
$(e_1, e_2, \ldots, e_i, \ldots, e_n)$.\cite{Brueckner1900}\\

\textbf{Example} $n=6$
\begin{figure}[!htp]
\begin{center}
\includegraphics[width=0.4\textwidth]{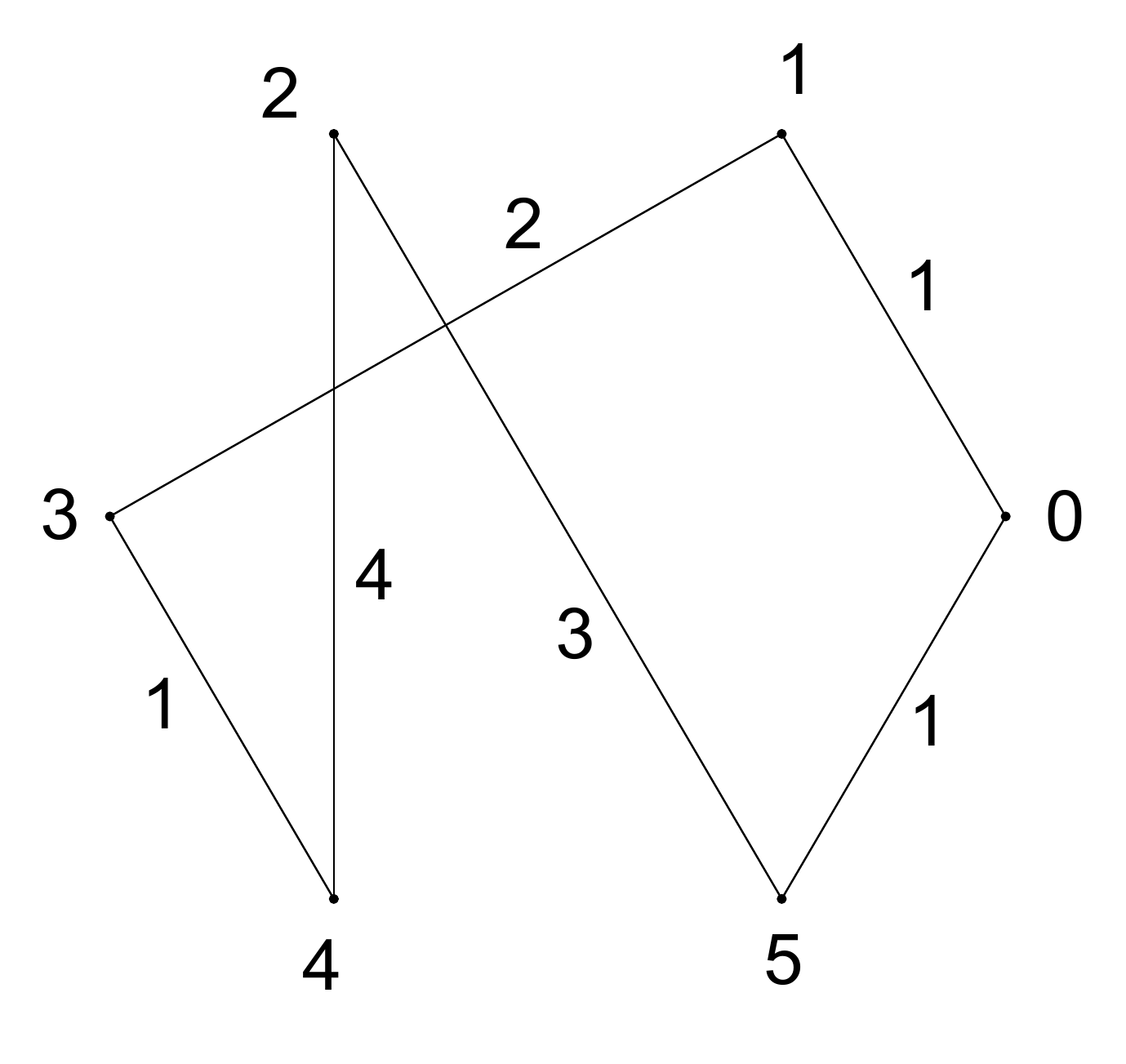}
\caption{Cycle of vertices (0 1 3 4 2 5) and cycle of edges (1 2 1 4 3 1)}
\end{center}
\end{figure}
\subsection{Definition of the basic equivalence relations and general formulas}
\label{subsec:definition_of_the_basic_equivalence_relations_and_general_formulas}
We denote by $C(n)$ the set of all $n$-polygons and define the following equivalence relations on $C(n)$:
\begin{enumerate}
\item[(1)] Two $n$-polygons $P_1(n)$ and $P_2(n)$ are said to be similar, denoted $P_1(n)\stackrel{\equiv}{_S}P_2(n)$, if they are obtainable from one another by a rotation or reflection.
\item[(2)] Two $n$-polygons $P_1(n)$ and $P_2(n)$ are said to be equivalent, denoted $P_1(n)\stackrel{\equiv}{_E}P_2(n)$, if they are obtainable from one another by a rotation, but not by a reflection.
\end{enumerate}
\newpage
\textbf{Example} $n=6$
\begin{figure}[!htp]
\begin{center}
\includegraphics[width=0.4\textwidth]{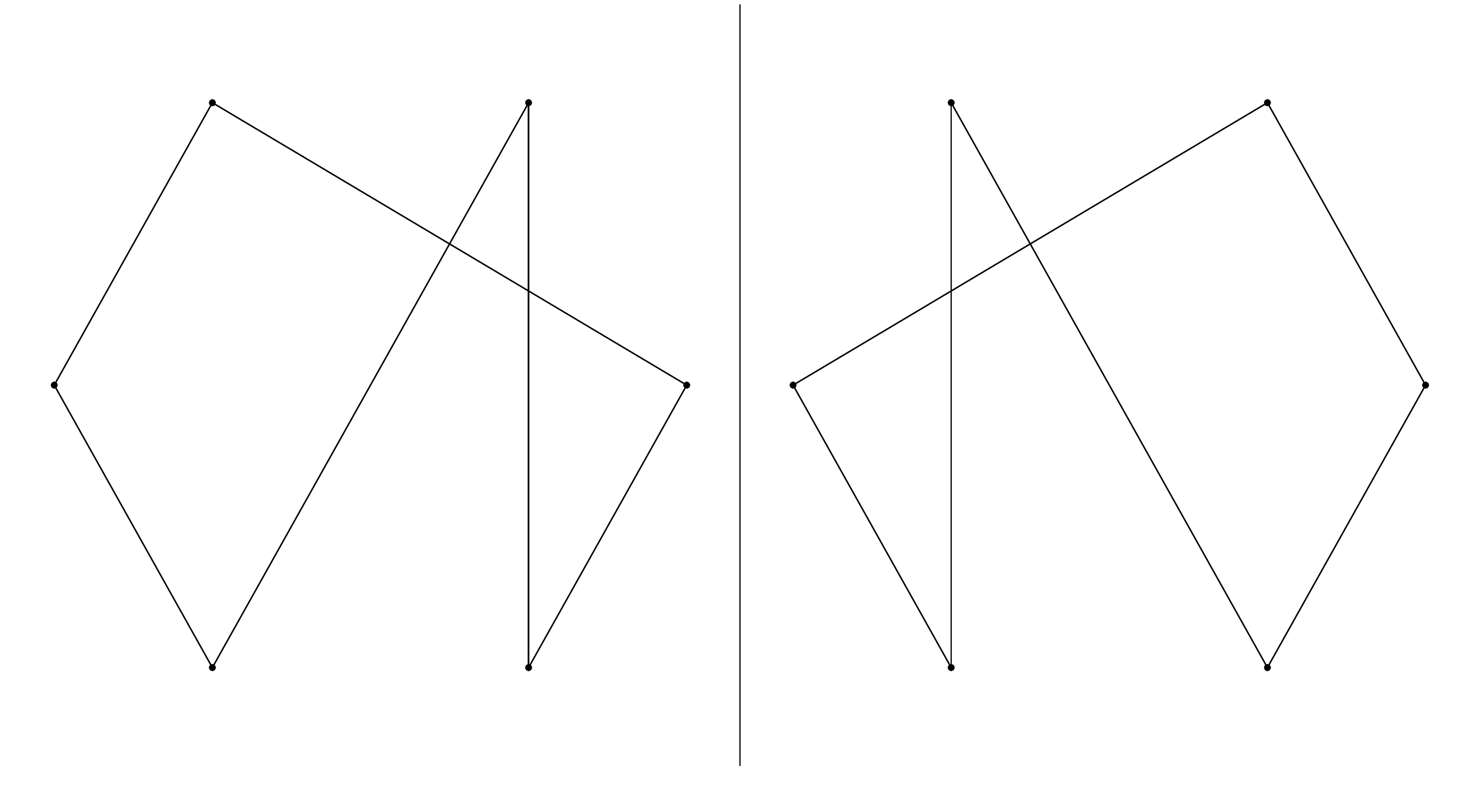}
\caption{The two $n$-polygons are similar but not equivalent.}
\end{center}
\end{figure}

The question of the number $\vert C(n)\stackrel{\equiv}{_S}\vert$ of equivalence classes of the similar $n$-polygons and the question of the number $\vert C(n)\stackrel{\equiv}{_E}\vert$ of equivalence classes of the equivalent $n$-polygons has been answered by Golomb/Welch \cite{Golomb1960} and proved by Herman/Poimenidou \cite{Herman2019} in a second way. So it is known how many different shapes of n-polygons exist for a given n.

\begin{enumerate}
\item $n$ odd: $\vert C(n)\stackrel{\equiv}{_E}\vert = \dfrac{1}{2n^2}\left(   \sum \limits_{d \mid n}  \varphi^2 \left( \dfrac{n}{d} \right)\cdot d! \cdot \left( \frac{n}{d} \right)^d \right)$
\item  $n$ odd: $\vert C(n)\stackrel{\equiv}{_S}\vert = \dfrac{1}{4n^2}\left(   \sum \limits_{d \mid n} \varphi^2 \left( \dfrac{n}{d} \right)\cdot d! \cdot \left( \frac{n}{d} \right)^d +2^{(n-1)/2} \cdot n^2 \cdot  \left( \frac{n-1}{2} \right)! \right)$
\item $n$ even: $\vert C(n)\stackrel{\equiv}{_E}\vert= \dfrac{1}{2n^2}\left(   \sum \limits_{d \mid n}  \varphi^2 \left( \dfrac{n}{d} \right)\cdot d! \cdot \left( \frac{n}{d} \right)^d + 2^{n/2} \cdot \left( \frac{n}{2} \right) \cdot \left( \frac{n}{2} \right)!\right)$
\item  $n$ even: $\vert C(n)\stackrel{\equiv}{_S}\vert = \dfrac{1}{4n^2}\left(   \sum \limits_{d \mid n} \varphi^2 \left( \dfrac{n}{d} \right)\cdot d! \cdot \left( \frac{n}{d} \right)^d +2^{n/2} \cdot \dfrac{n(n+6)}{4} \cdot  \left( \frac{n}{2} \right)! \right)$
\end{enumerate}

In these formulas, $\varphi \left(\dfrac{n}{d} \right) $ denotes the Euler $\varphi$ function of $\frac{n}{d}$ and $d$ a divisor of $n$.\\

\subsection{The question to deal with in this article}
\label{subsec:The_question_to_deal_with_in_this_article}

A closer look at the different shapes reveals that they have different and only certain symmetry properties: E.g. for $n = 6$:
\begin{enumerate}
\item[$\bullet$] One single shape with six axes,
\item[$\bullet$] one shape with three axes,
\item[$\bullet$] three different shapes with two axes,
\item[$\bullet$] five different shapes with one axis,
\item[$\bullet$] four different shapes without any axis.
\end{enumerate}
 
In the latter case we observe that two shapes are completely asymmetrical, but they are similar. Two other ones can be made to coincide to themselves by turning them 180 degrees around the center of the circle. They are also similar.
\newpage

\textbf{Example} $n=6$
\begin{figure}[!h]
\begin{center}
\begin{tabular}{c | c | c | c}
\includegraphics[width=0.2\textwidth]{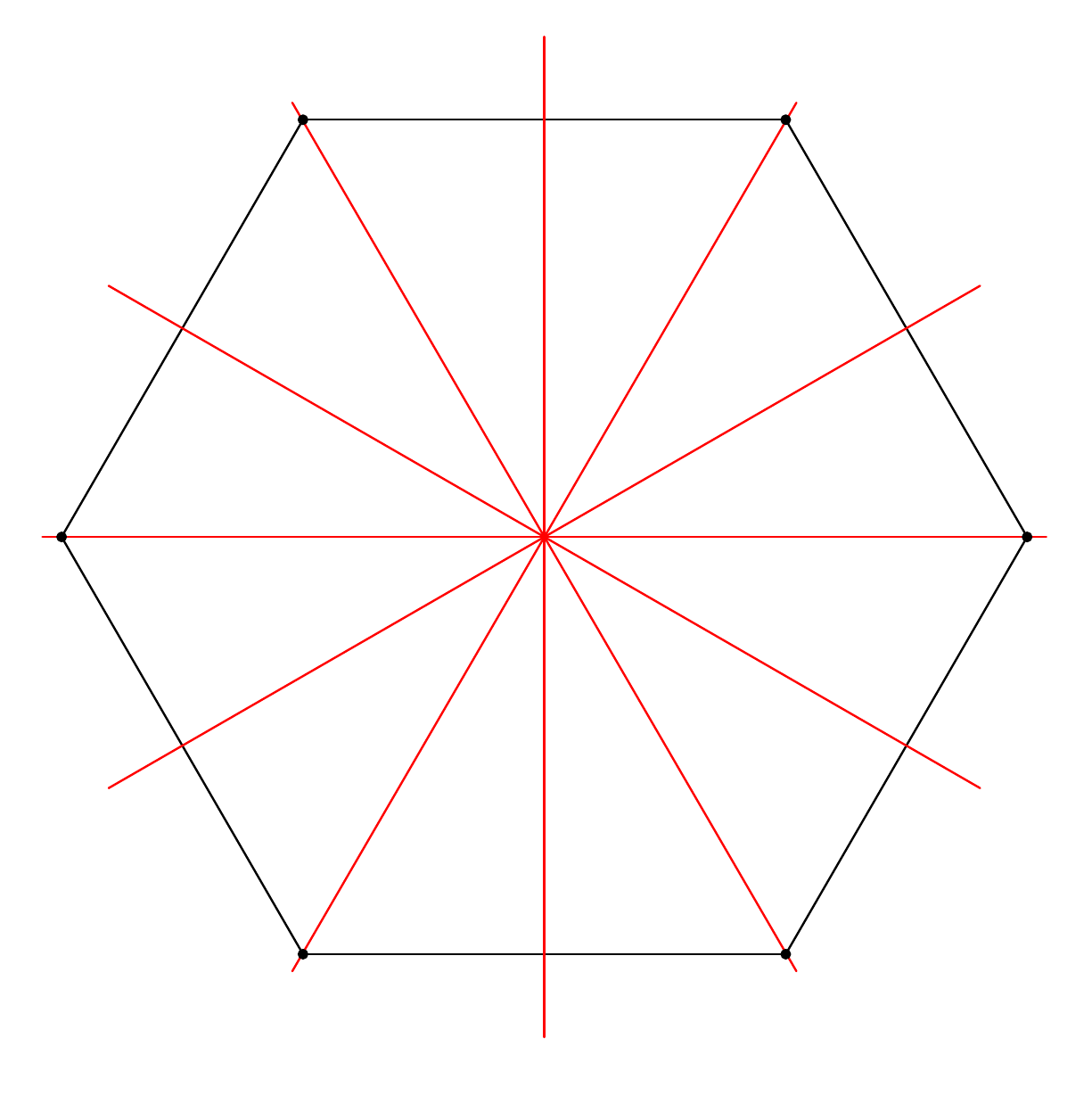} & \includegraphics[width=0.2\textwidth]{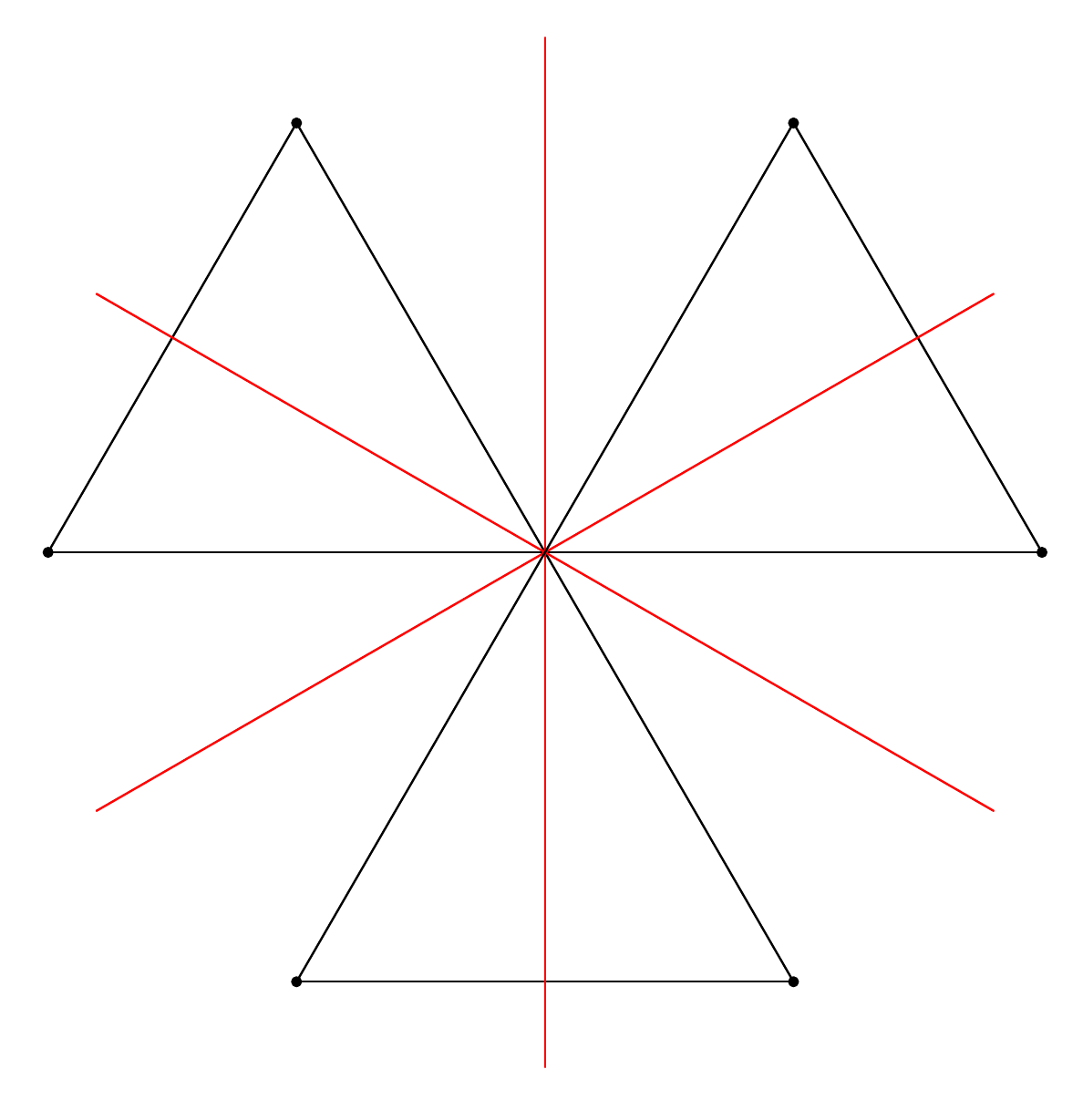} & \includegraphics[width=0.2\textwidth]{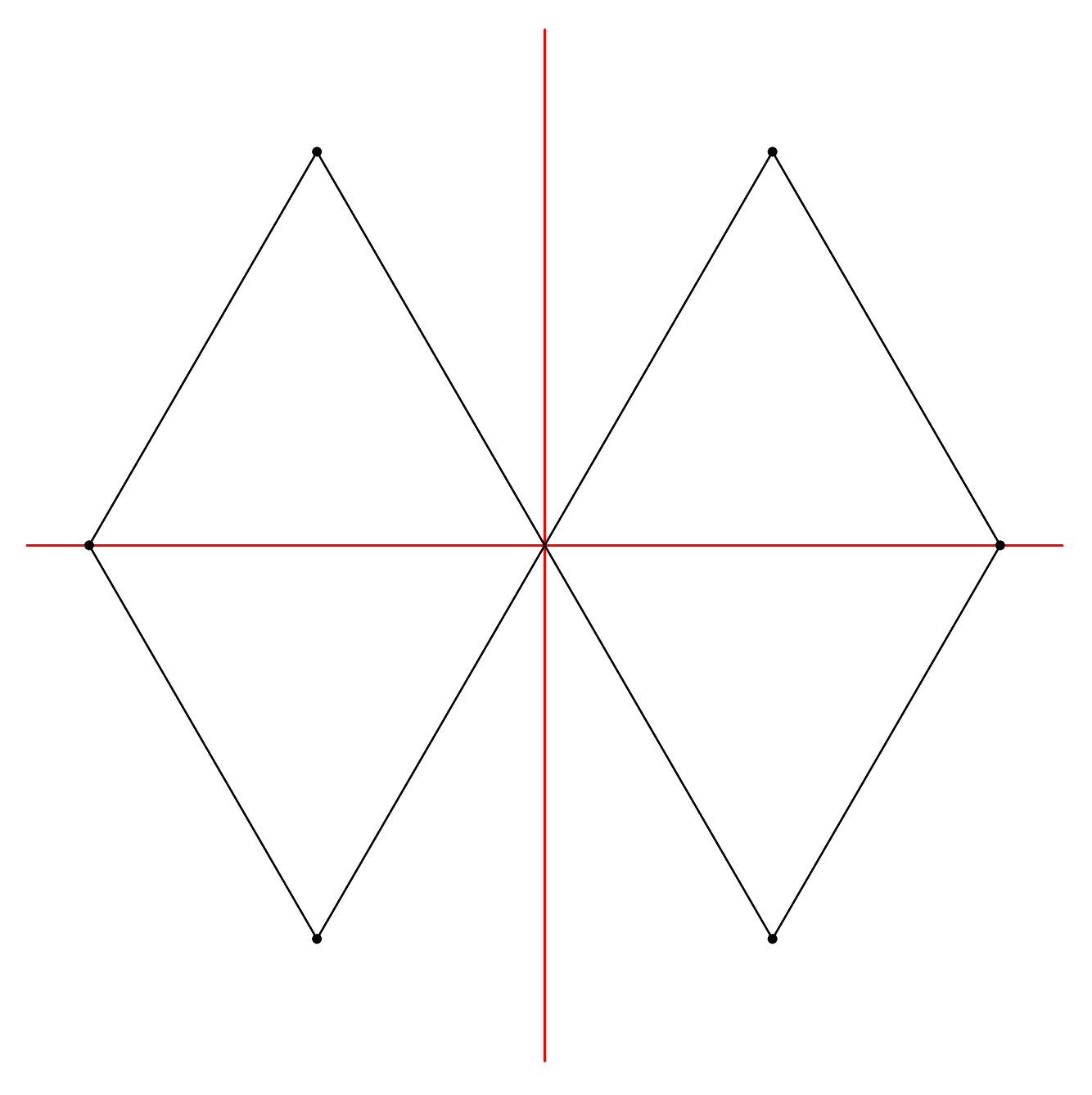} & \includegraphics[width=0.2\textwidth]{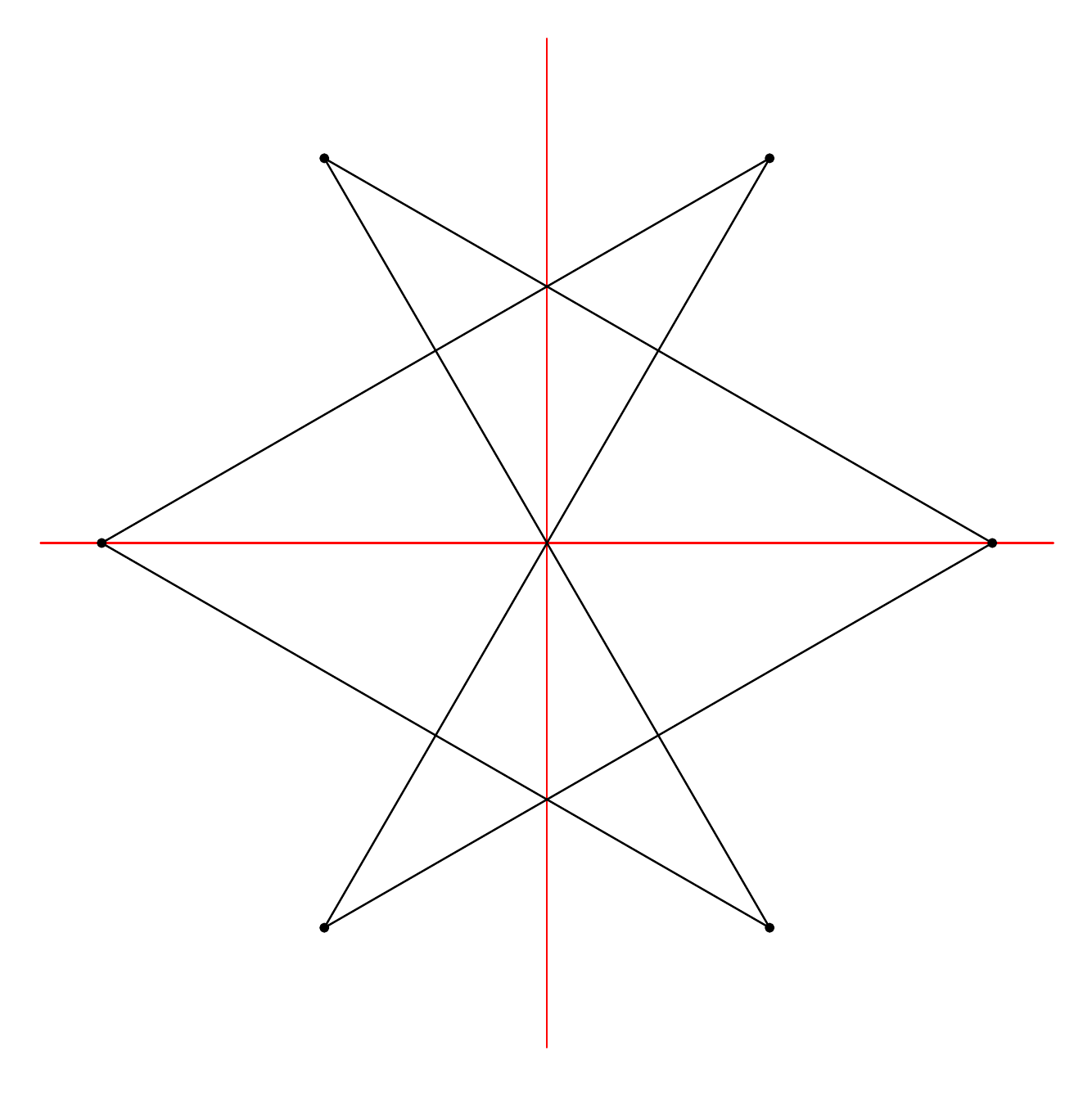}\\
6 axes & 3 axes & 2 axes & 2 axes\\ \hline 
\includegraphics[width=0.2\textwidth]{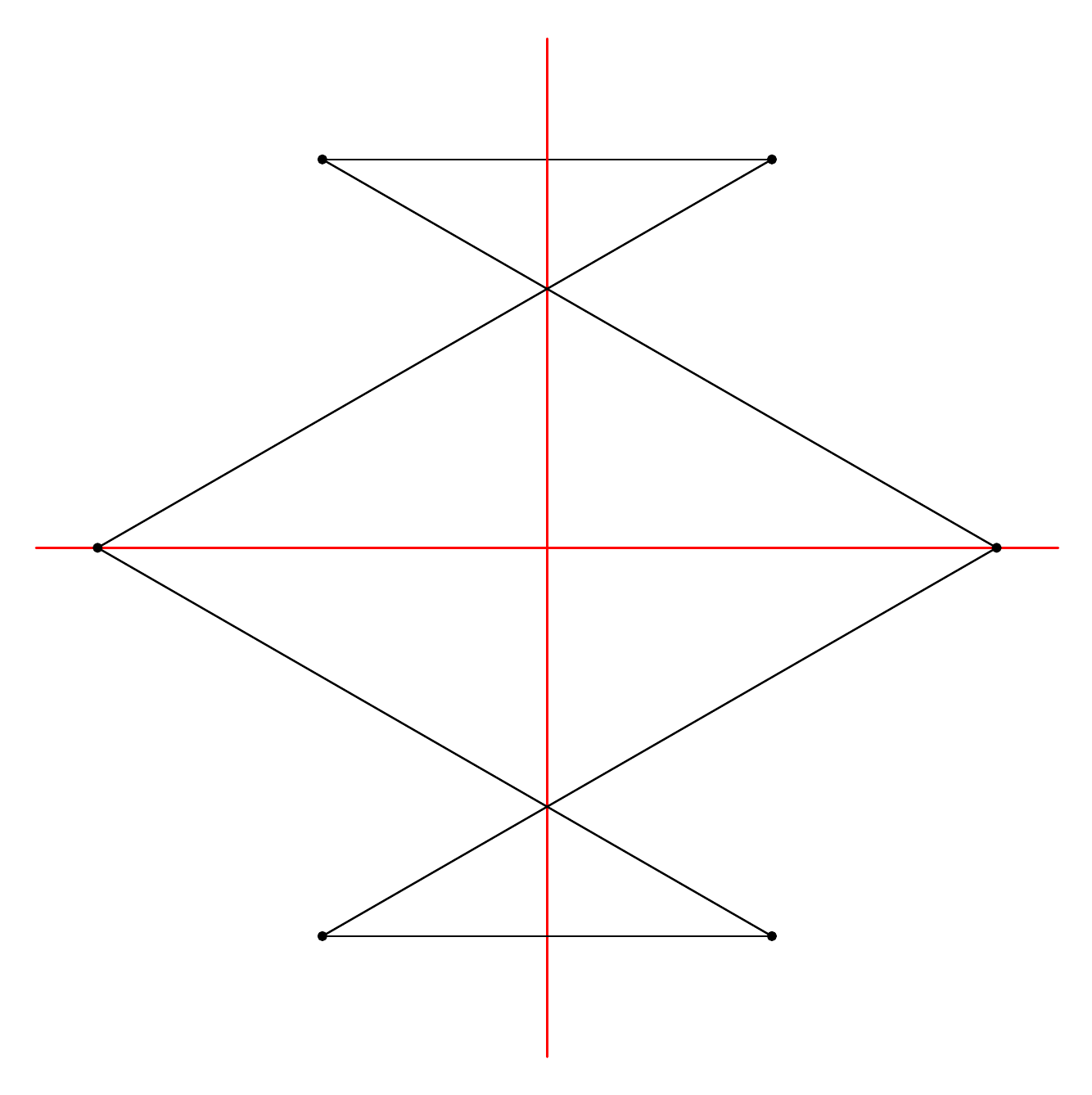} & \includegraphics[width=0.2\textwidth]{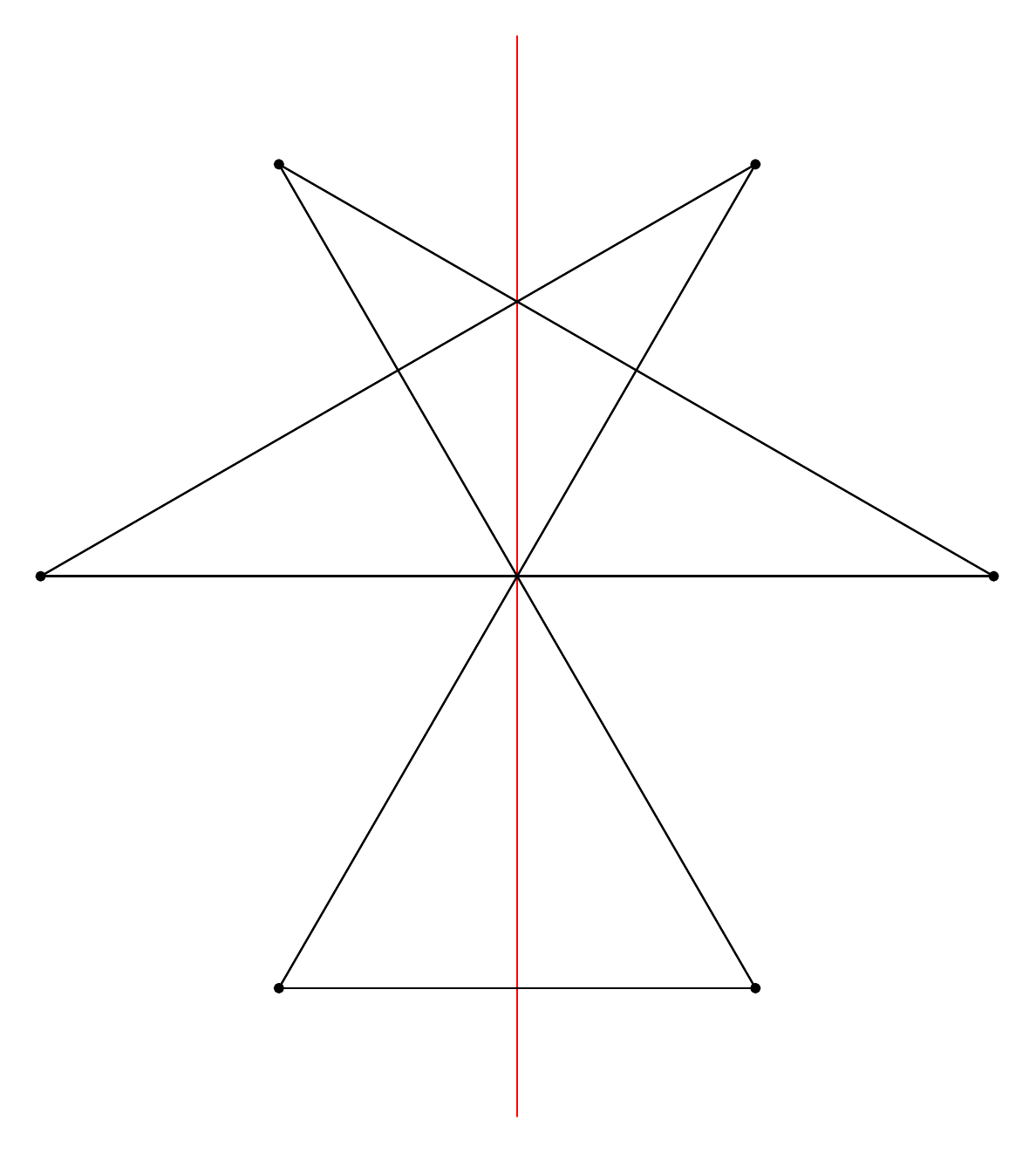} & \includegraphics[width=0.2\textwidth]{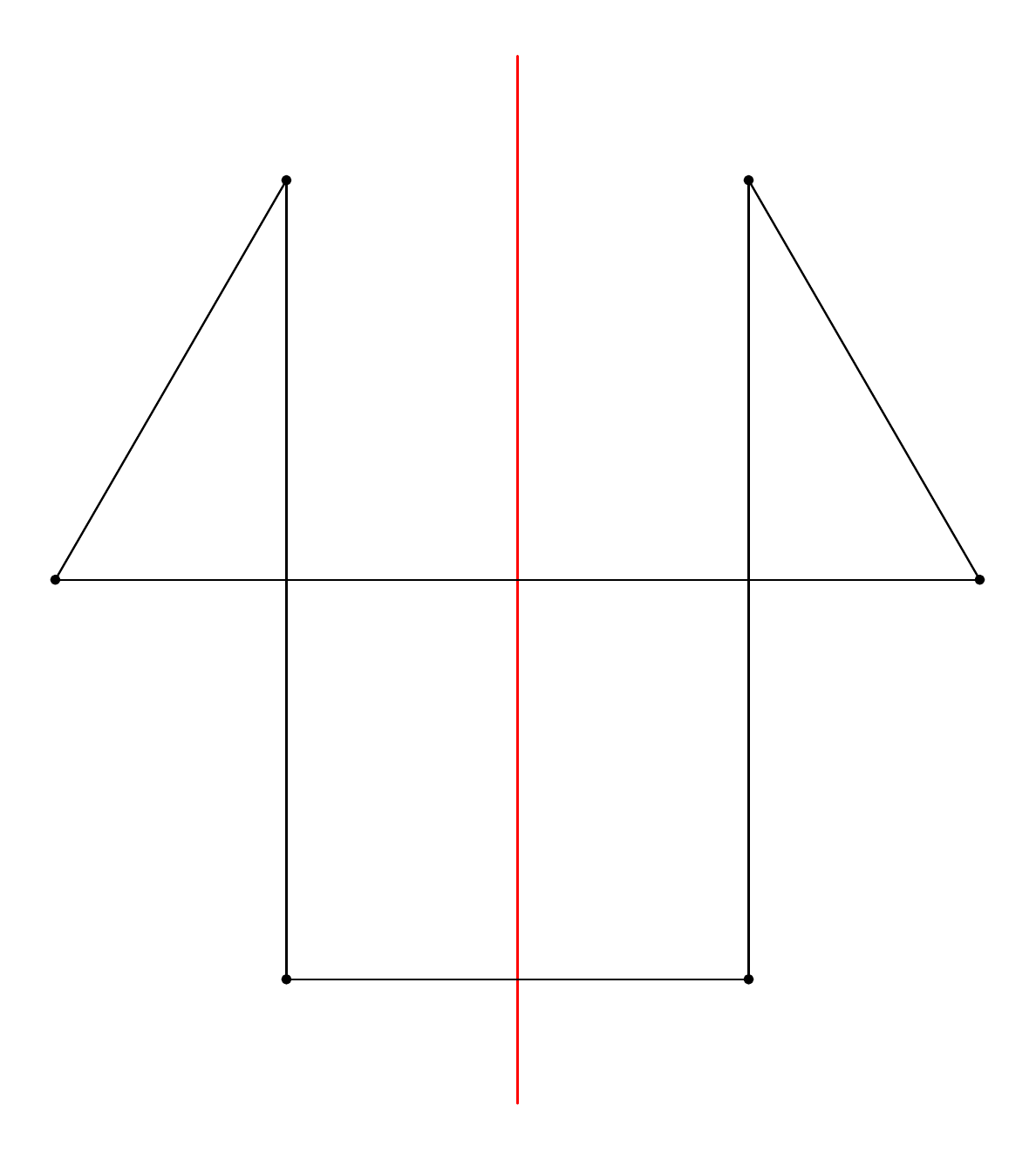} & \includegraphics[width=0.2\textwidth]{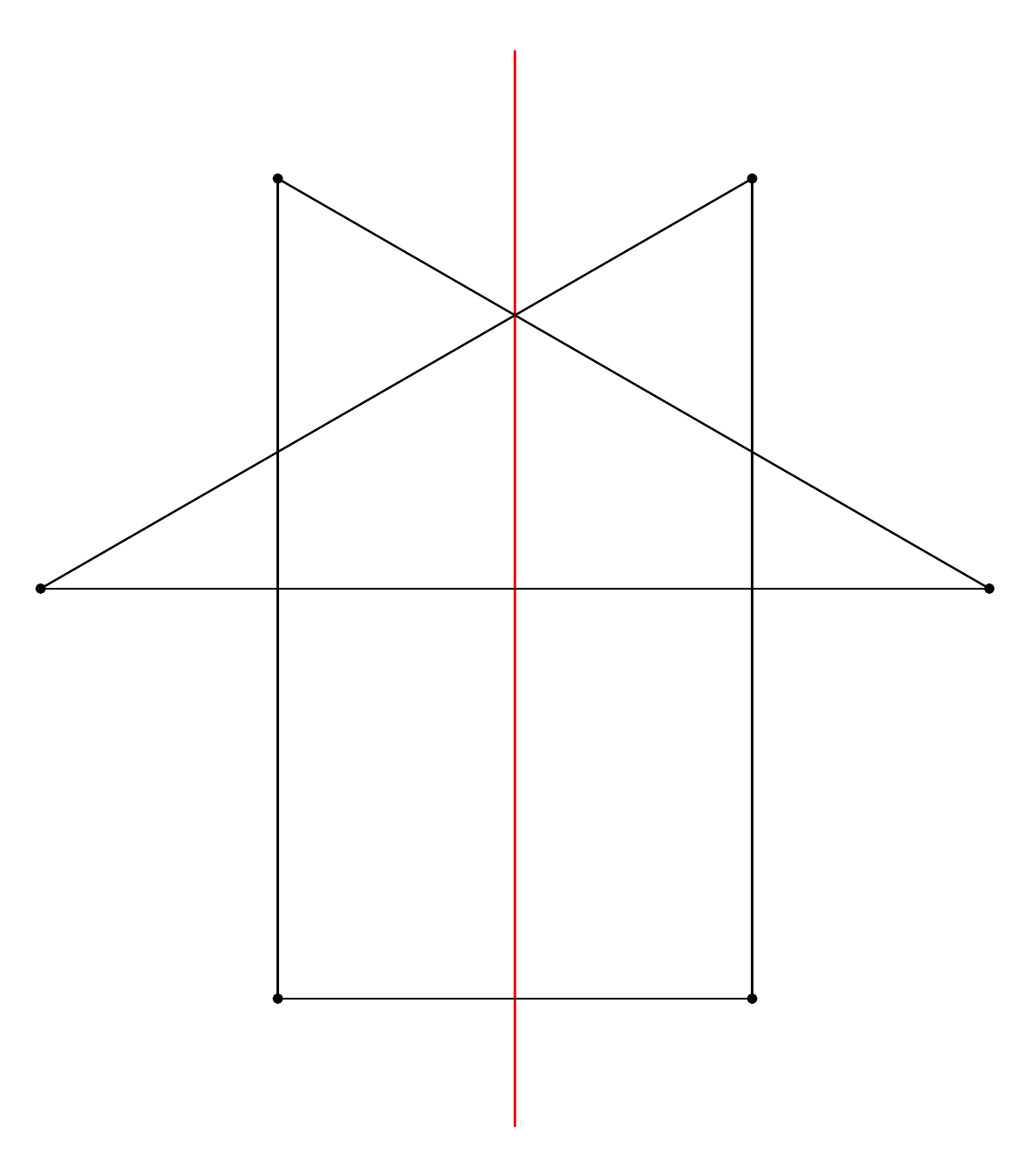}\\
2 axes & 1 axis & 1 axis & 1 axis\\ \hline 
\includegraphics[width=0.2\textwidth]{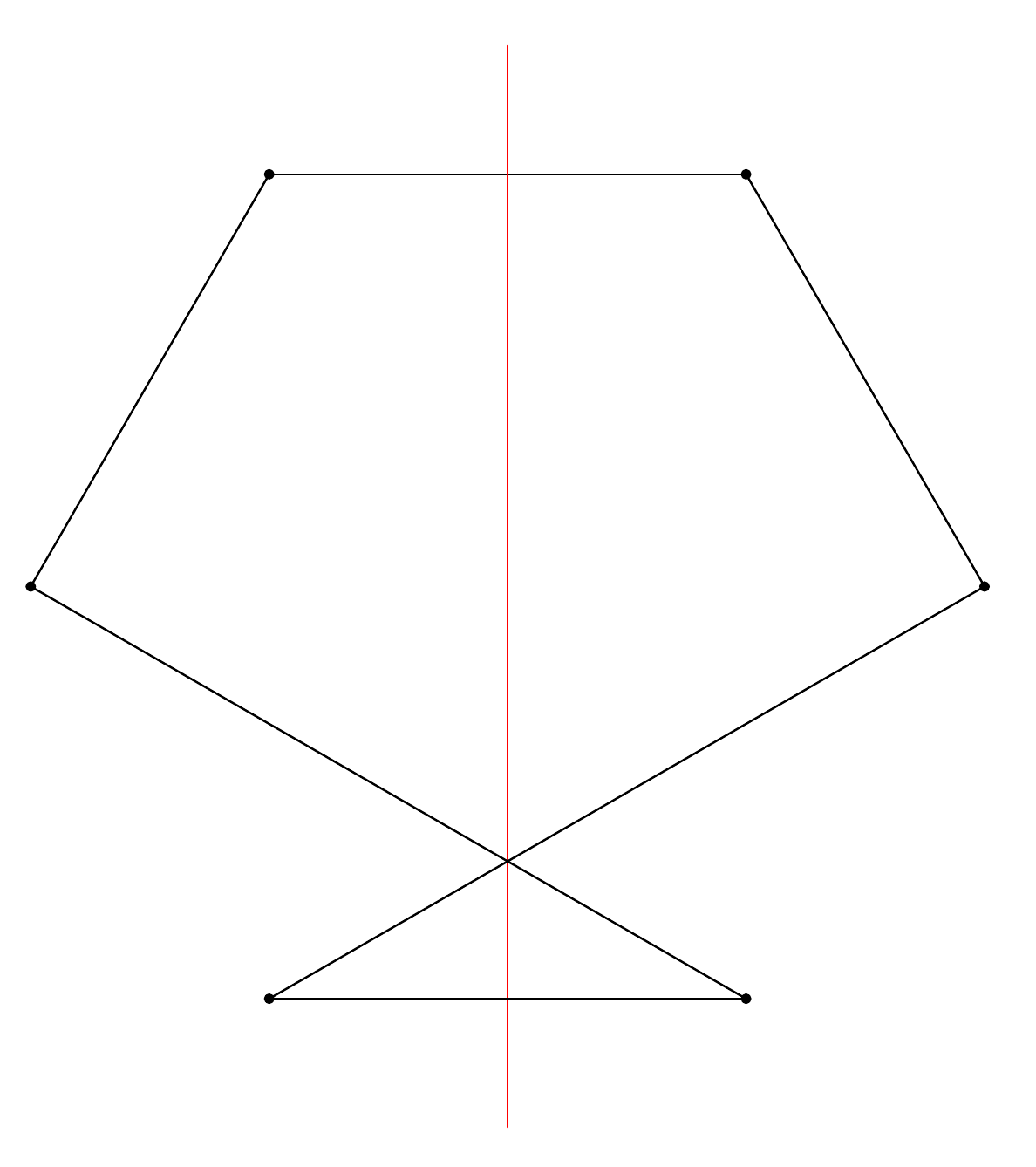} & \includegraphics[width=0.2\textwidth]{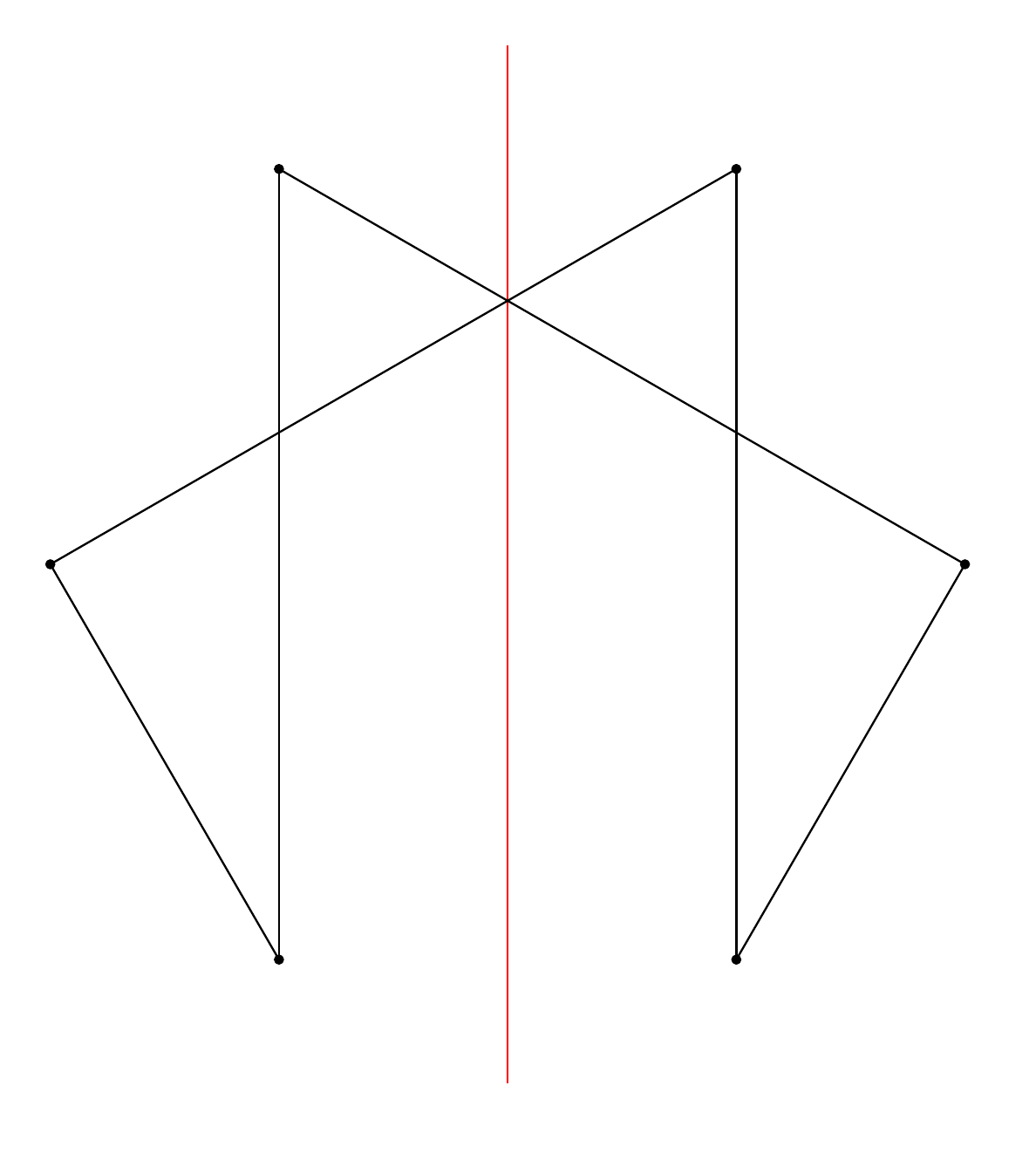} &  & \\
1 axis & 1 axis &  & \\ \hline
\includegraphics[width=0.2\textwidth]{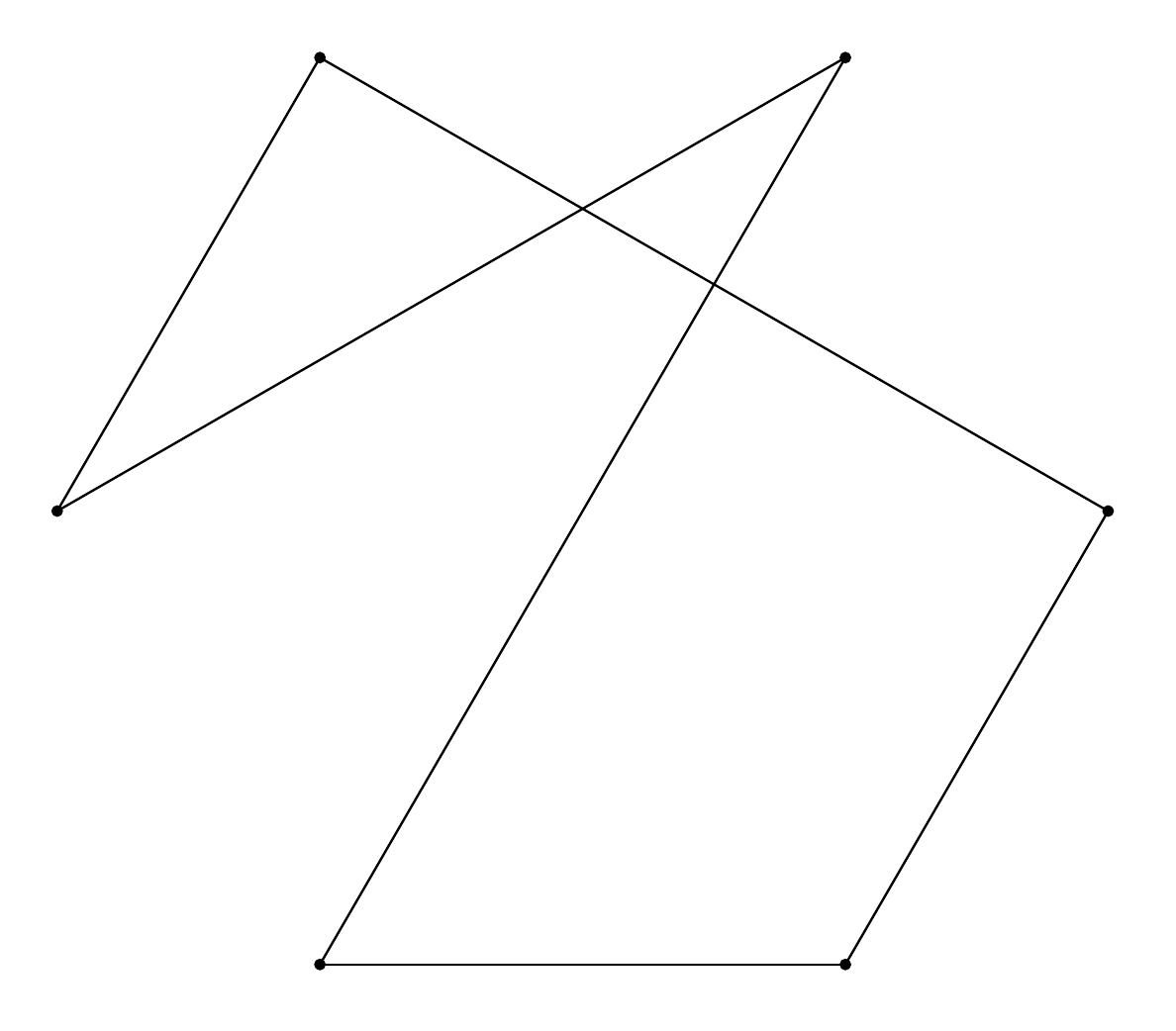} & \includegraphics[width=0.2\textwidth]{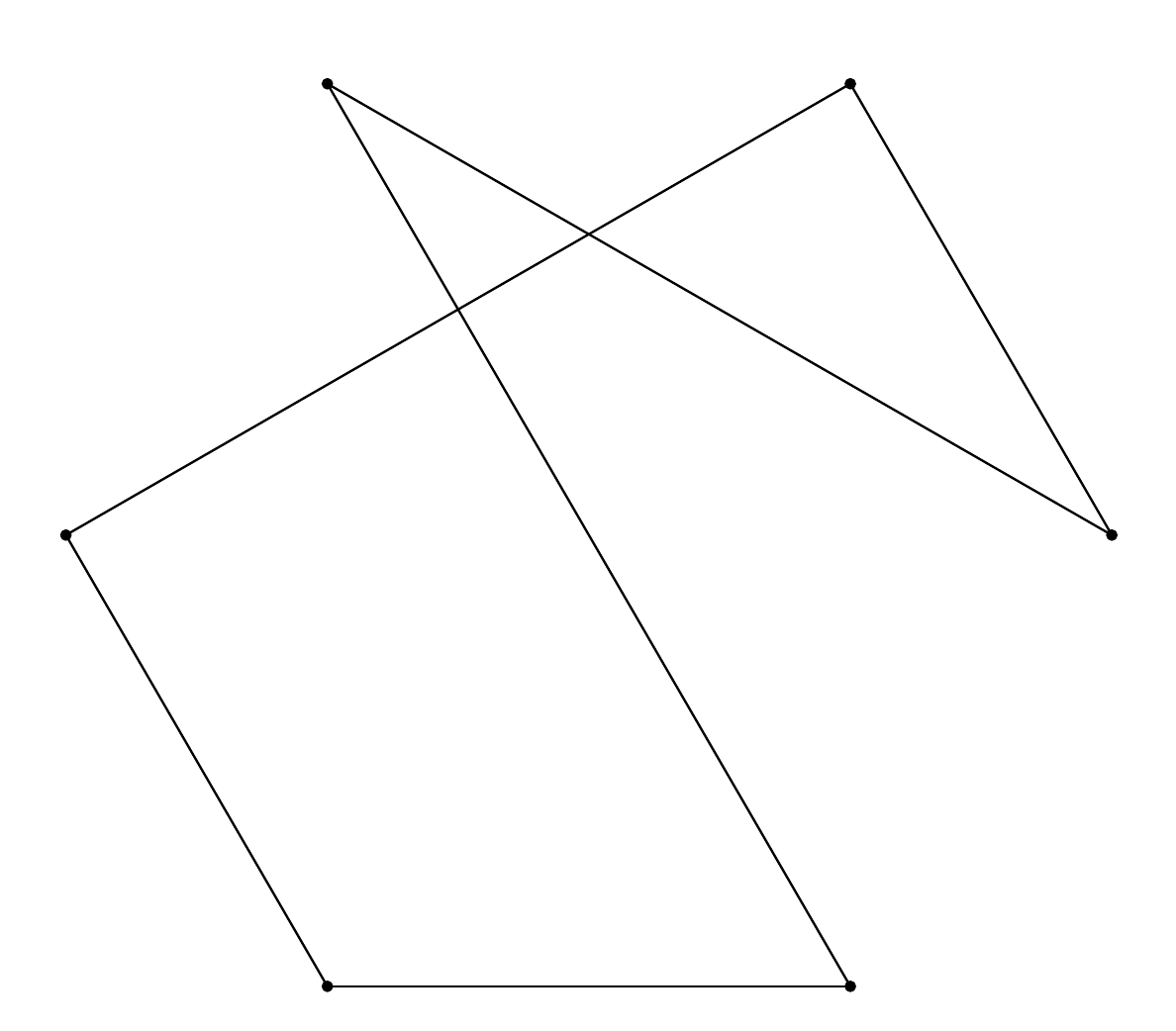} & \includegraphics[width=0.2\textwidth]{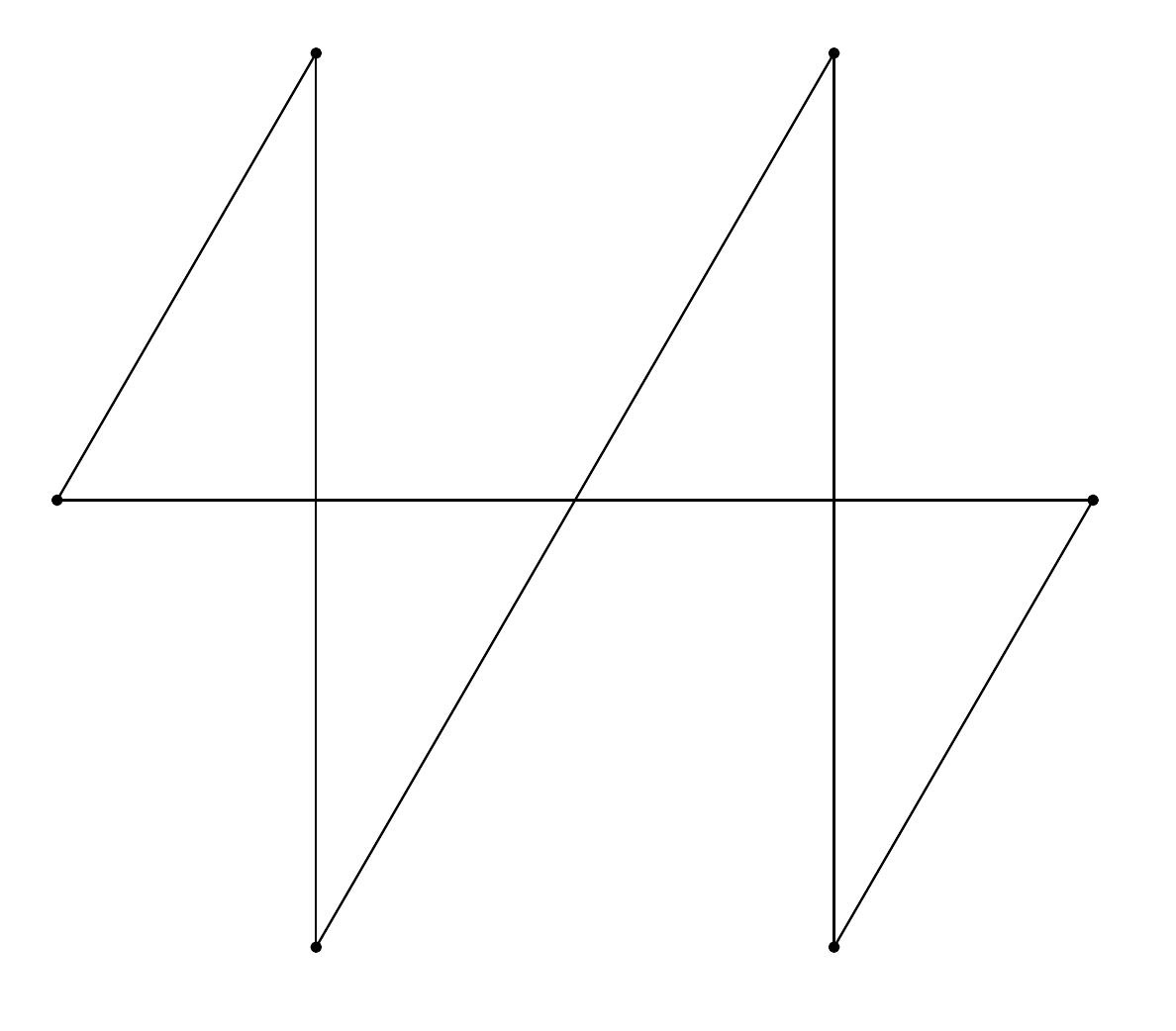} & \includegraphics[width=0.2\textwidth]{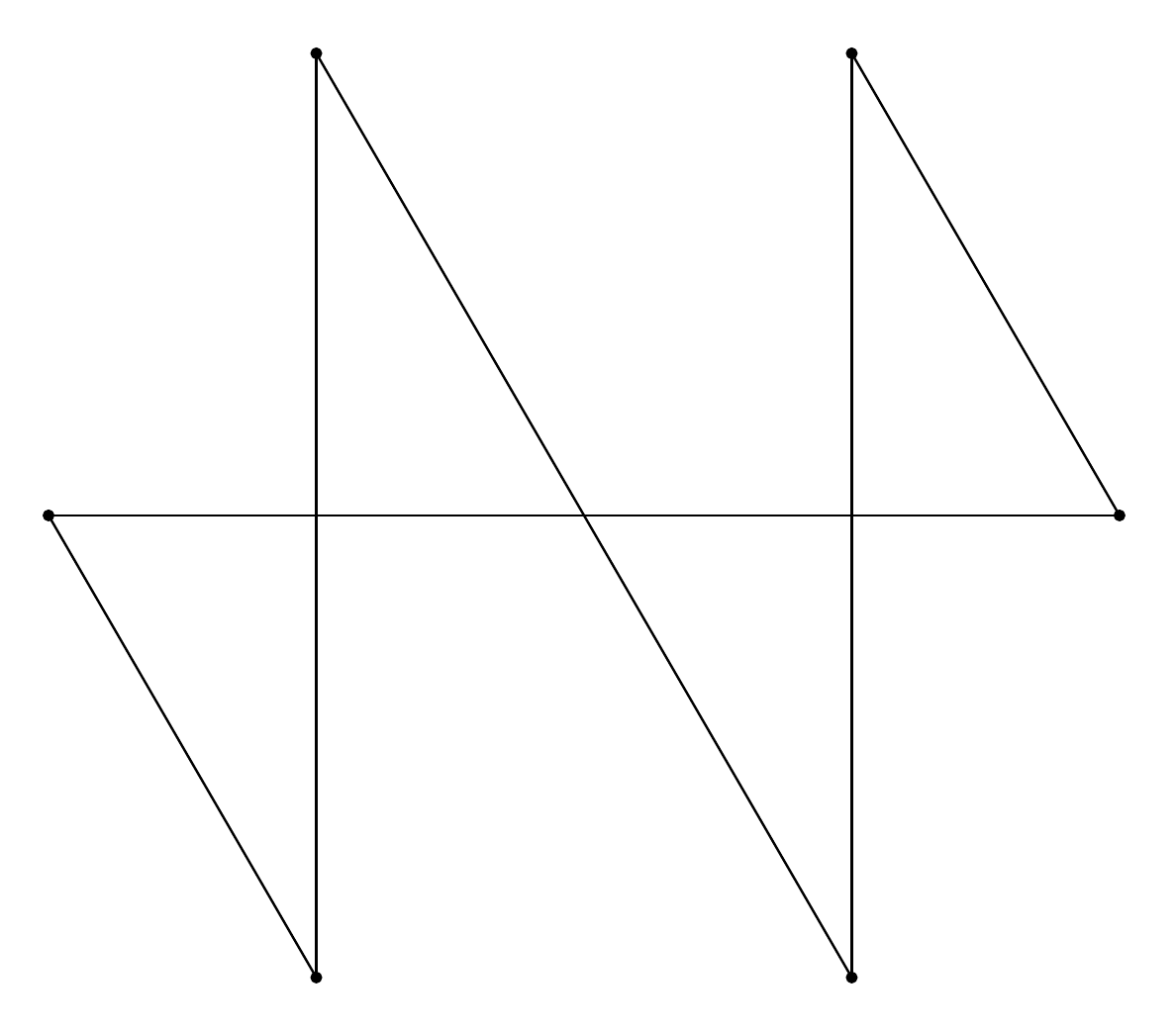}\\
no axes & no axes & no axes & no axes\\ \hline
\end{tabular}
\caption{A set of representatives for $n=6$, sorted by symmetry-level}
\end{center}
\end{figure}

Such observations suggest the following definition and question:\\

\textbf{Definition of a $m$-symmetric $n$-polygon} 

Let $m>0$ be an integer. A $m$-symmetric $n$-polygon is a $n$-polygon with $m$ axes, denoted by $P_m(n)$. $\vert P_m(n) \vert$ denotes the number of their equivalence classes.\\

\textbf{Question of the article}

Let $n = 2m$ be an even number. We ask for the number $\vert P_m(n) \vert$ of equivalence classes of the $m$-symmetric $n$-polygons.

\section{Results}
\label{sec:results}
\subsection{Main theorem}
\label{subsec:main_theorem}

Let $n=2m>3$ be an even integer.\newline
The different equivalence classes of $n$-polygons with $m$ axes are represented by the $n$-tuples $(a, b, a, b, \ldots, a, b)$ of their sides, if and only if $a$ and $b$ have the following four properties:
\begin{enumerate}
\item $a\in \mathbb{N}$ with $a\equiv $1 mod 2,
\item $b\in \mathbb{N}$ with $b\equiv $1 mod 2,
\item $1 \leq a<b \leq n-1$,
\item $gcd\left(\dfrac{a+b}{2}=:u,m \right)=1$.
\end{enumerate}

\subsection{Conclusion: Formula for $\vert P_m(n) \vert$}
\label{subsec:conclusion_formula_for_p_m_n}

Let $\dfrac{a+b}{2}:= u$ be prime to $m$. For each allowed $u$-value, we determine all pairs $(a, b)$ which satisfy the first three properties in the main theorem. It turns out that the number of allowed pairs $(a, b)$ belonging to a certain $u$-value is very easily determinable: If $u$ is even, then the number of allowed pairs is $\dfrac{u}{2}$, if $u$ is odd, then the number of the allowed pairs is $\dfrac{u-1}{2}$. This yields to the following final formula for the number $\vert P_m(n) \vert$ of equivalence classes of $n$-polygons with $m$ axes:

\begin{center}
$\vert P_m(n) \vert=\sum\limits_{\substack{u\equiv 0 mod 2,\\gcd \left( u,m \right)=1 }} \frac{u}{2}+\sum\limits_{\substack{u\equiv 1 mod 2,\\gcd \left( u,m \right)=1 }} \frac{u-1}{2}$
\end{center}

\begin{table}[!htp]
\centering
\begin{tabular}{| c | c | c || c | c | c || c | c | c || c | c | c |}
\hline
n & m & $\vert P_m(n) \vert$ & n & m & $\vert P_m(n) \vert$ & n & m & $\vert P_m(n) \vert$ & n & m & $\vert P_m(n) \vert$\\ \hline
4 & 2 & 1 & 6 & 3 & 1 & 8 & 4 & 1 & 10 & 5 & 4 \\
12 & 6 & 2 & 14 & 7 & 9 & 16 & 8 & 6 & 18 & 9 & 12\\
20 & 10 & 8 & 22 & 11 & 25 & 24 & 12 & 10 & 26 & 13 & 36\\
28 & 14 & 18 & 30 & 15 & 28 & 32 & 16 & 28 & 34 & 17 & 64\\
36 & 18 & 24 & 38 & 19 & 81 & 40 & 20 & 36 & 42 & 21 & 60\\
44 & 22 & 50 & 46 & 23 & 121 & 48 & 24 & 44 & 50 & 25 & 120\\
52 & 26 & 72 & 54 & 27 & 117 & 56 & 28 & 78 & 58 & 29 & 196\\
60 & 30 & 56 & 62 & 31 & 225 & 64 & 32 & 120 & 66 & 33 & 160\\
68 & 34 & 128 & 70 & 35 & 204 & 72 & 36 & 102 & 74 & 37 & 324\\
76 & 38 & 162 & 78 & 39 & 228 & 80 & 40 & 152 & 82 & 41 & 400\\
84 & 42 & 120 & 86 & 43 & 441 & 88 & 44 & 210 & 90 & 45 & 264\\
\hline
\end{tabular}
\caption{Number of equivalence classes of $n$-polygons with $m$ axes}
\label{tab:number_of_equivalence_classes_of_n-polygons_with_m_axes}
\end{table}

\subsection{Example: $n=30$}
\label{subsec:n_30}
There are $\varphi(15)=15 \cdot \left(1-\dfrac{1}{3} \right)\cdot\left(1-\dfrac{1}{5} \right)=8$ numbers, which are prime to 15. They are 1, 2, 4, 7, 8, 11, 13 and 14. To each of these $u$-values we determine the number of the pairs $(a,b)$, which satisfy the first three properties of the main theorem.\\
\newpage
\begin{table}[!h]
\centering
\begin{tabular}{| c | c | c || c | c | c || c | c | c || c | c | c |}
\hline
u & a & b & u & a & b & u & a & b & u & a & b\\ \hline
1 & - & - & 2 & 1 & 3 & 4 & 1 & 7 & 7 & 1 & 13\\
- & - & - & - & - & - & 4 & 3 & 5 & 7 & 3 & 11\\
- & - & - & - & - & - & - & - & - & 7 & 5 & 9\\ \hline
- & 0 & - & - & 1 & - & - & 2 & - & - & 3 & -\\ \hline
\hline
u & a & b & u & a & b & u & a & b & u & a & b\\ \hline
8 & 1 & 15 & 11 & 1 & 21 & 13 & 1 & 25 & 14 & 1 & 27\\
8 & 3 & 13 & 11 & 3 & 19 & 13 & 3 & 23 & 14 & 3 & 25\\
8 & 5 & 11 & 11 & 5 & 17 & 13 & 5 & 21 & 14 & 5 & 23\\
8 & 7 & 9 & 11 & 7 & 15 & 13 & 7 & 19 & 14 & 7 & 21\\
- & - & - & 11 & 9 & 13 & 13 & 9 & 17 & 14 & 9 & 19\\
- & - & - & - & - & - & 13 & 11 & 15 & 14 & 11 & 17\\
- & - & - & - & - & - & - & - & - & 14 & 13 & 15\\ \hline
- & 4 & - & - & 5 & - & - & 6 & - & - & 7 & -\\ \hline
\end{tabular}\\
\caption{Allowed pairs $(a,b)$}
\label{tab:allowed_pairs }
\end{table}

It follows:
\begin{center}
$\vert P_{15}(30) \vert=0+1+2+3+4+5+6+7=\underline{\underline{28}}$
\end{center}

\subsection{A set of 28 representatives of the 30-polygons with 15 axes}
\label{subsec:a_set_of_28_representatives_of_the_30-polygons_with_15_axes}

\begin{figure}[!h]
\centering
\begin{tabular}{c | c | c | c}
\includegraphics[width=0.2\textwidth]{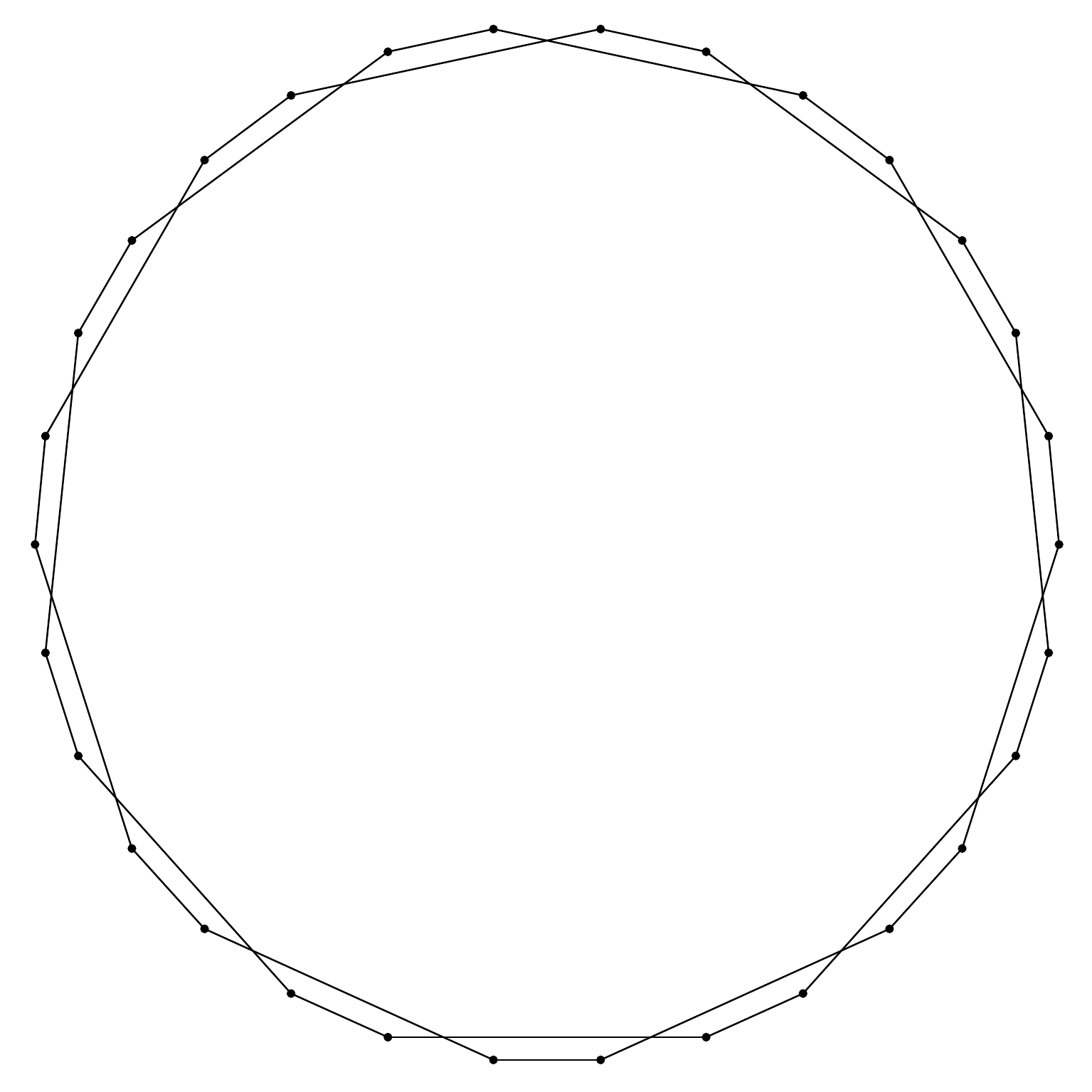} & \includegraphics[width=0.2\textwidth]{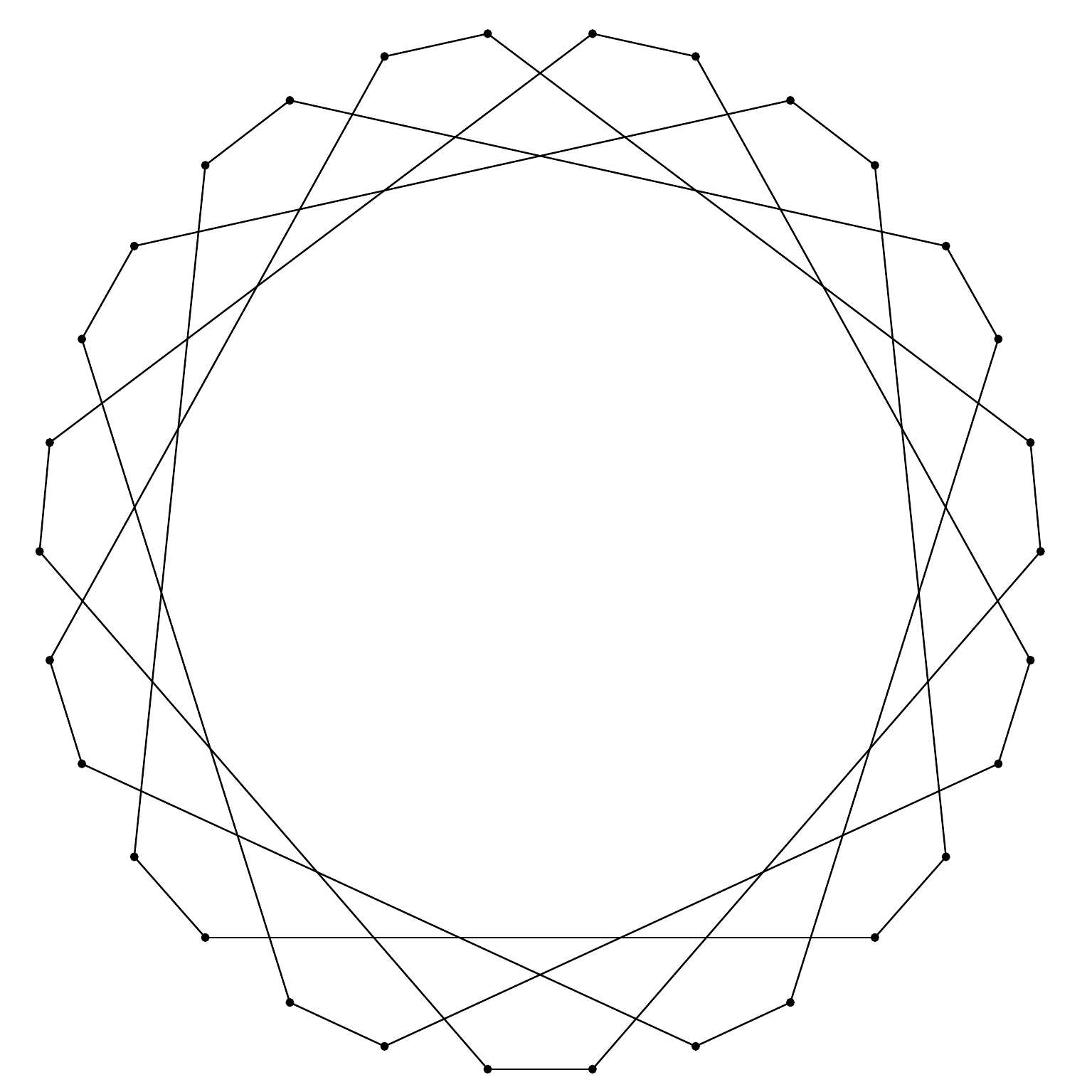} & \includegraphics[width=0.2\textwidth]{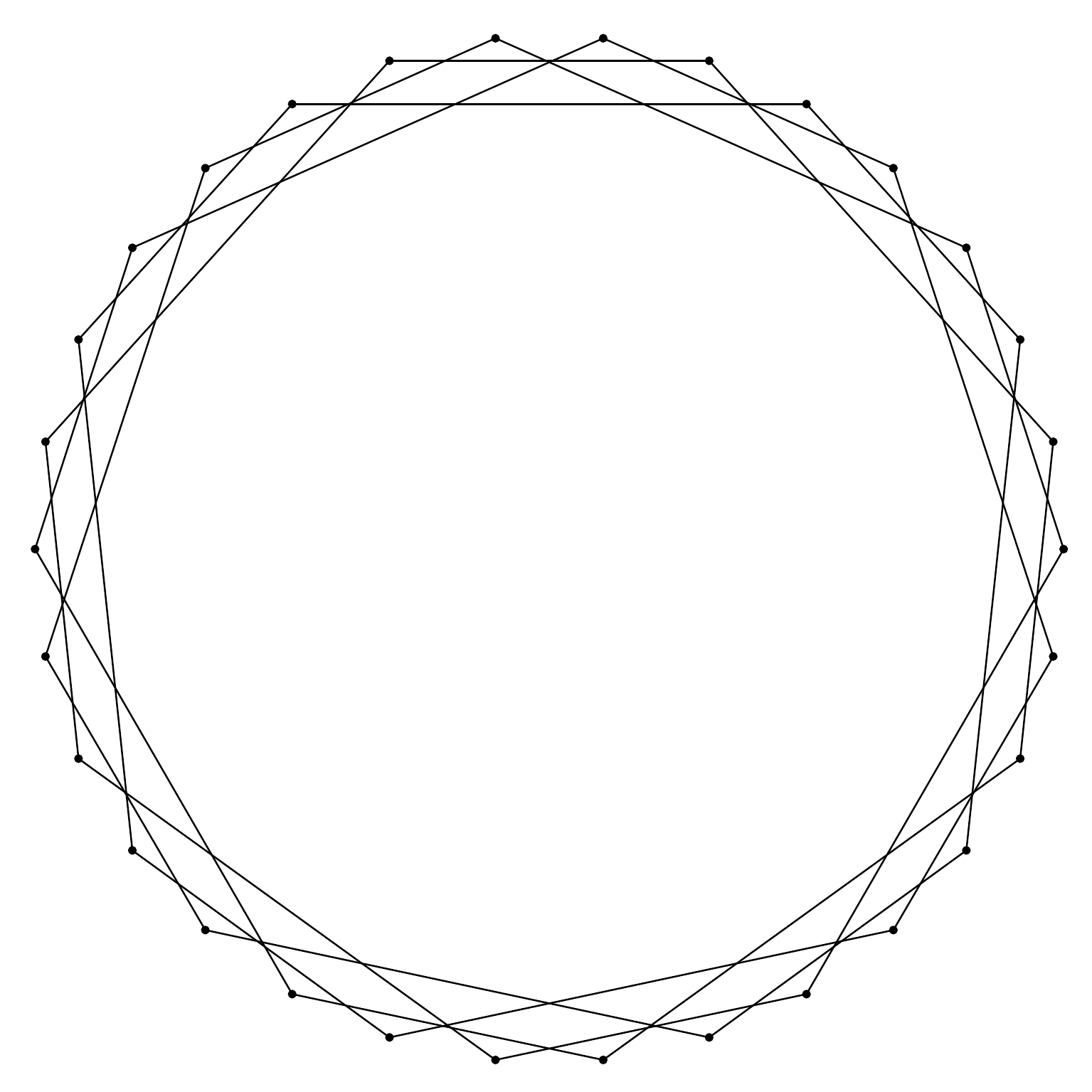} & \includegraphics[width=0.2\textwidth]{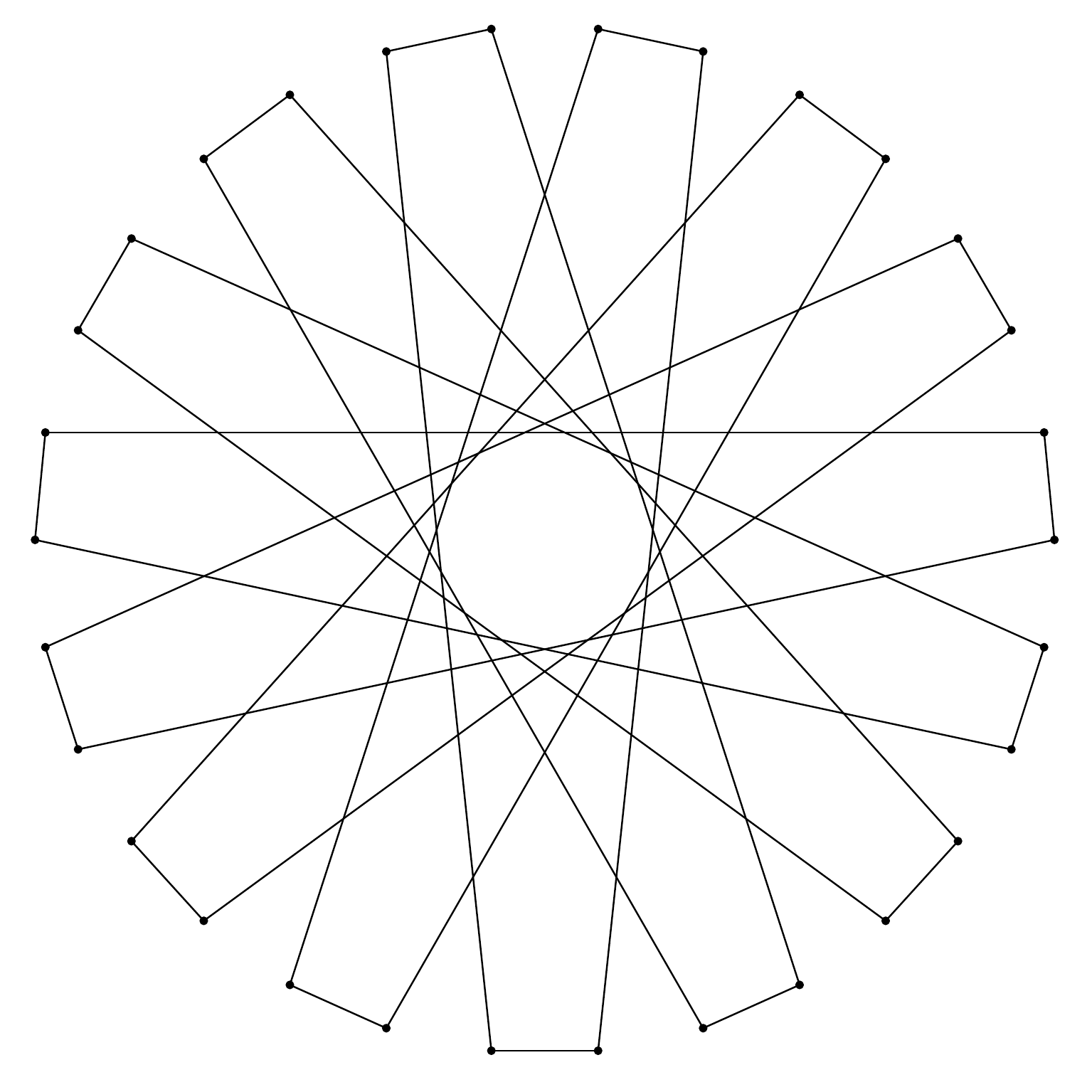}\\
$u=2;a=1;b=3$ & $u=4;a=1;b=7$ & $u=4;a=3;b=5$ & $u=7;a=1;b=13$\\ \hline
\includegraphics[width=0.2\textwidth]{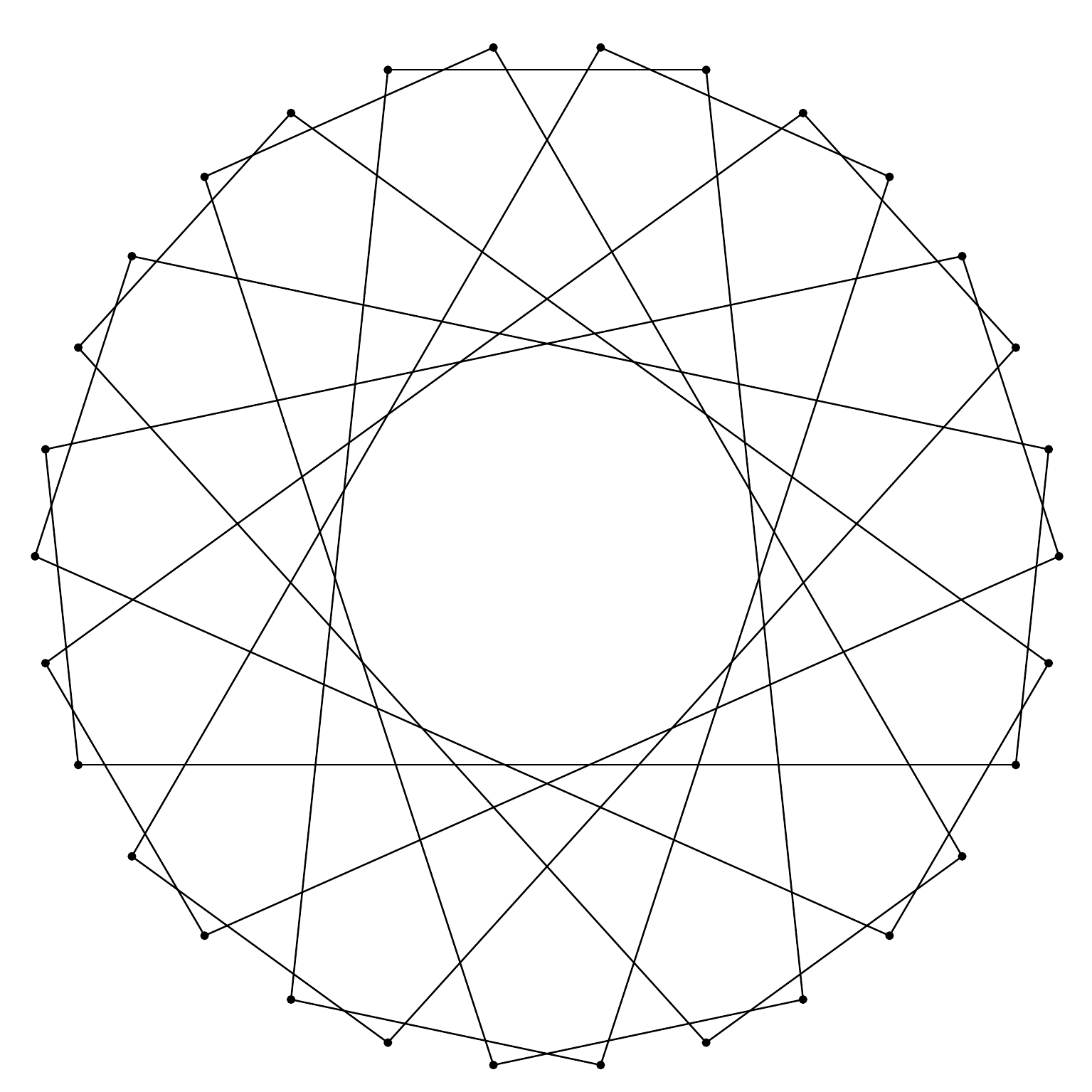} & \includegraphics[width=0.2\textwidth]{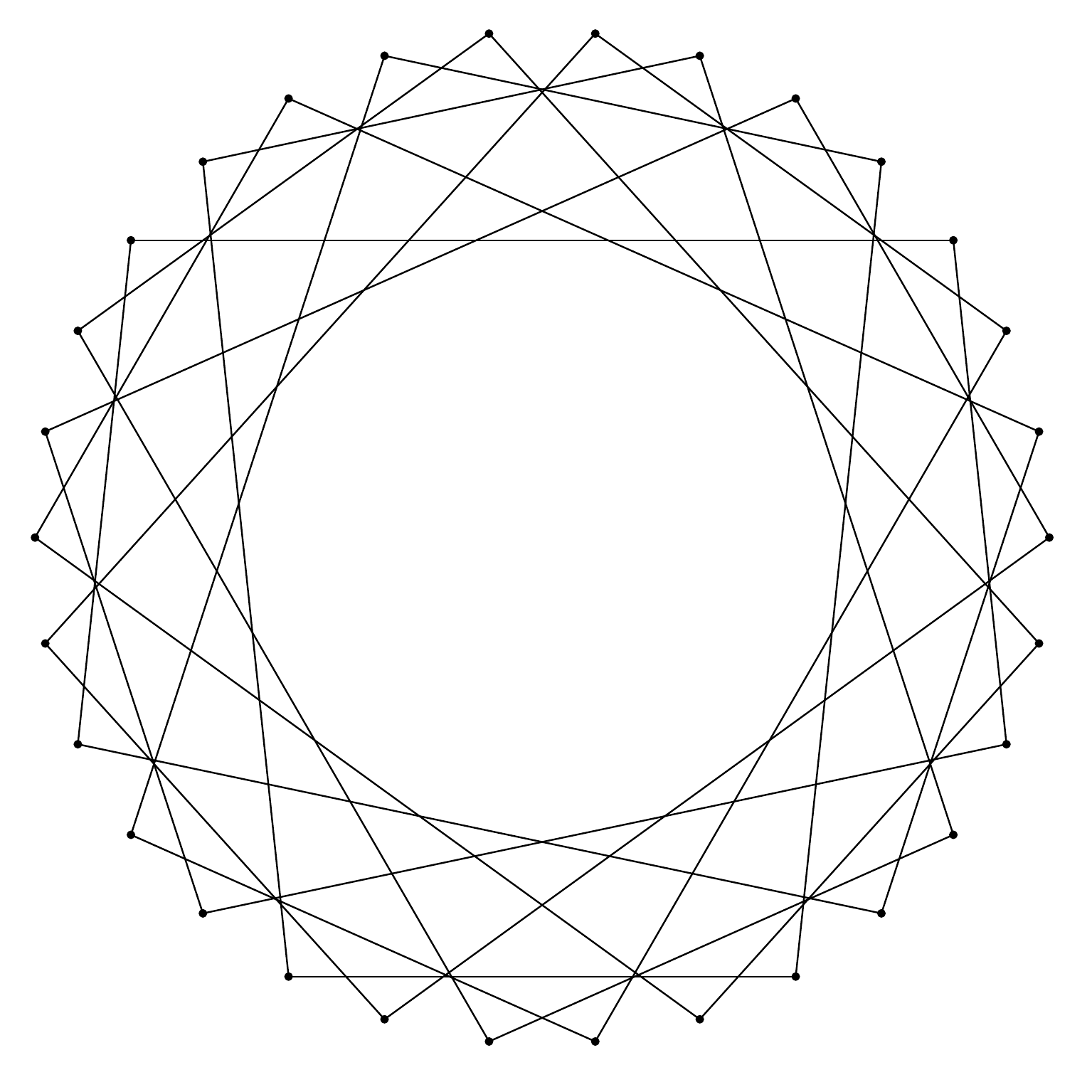} & \includegraphics[width=0.2\textwidth]{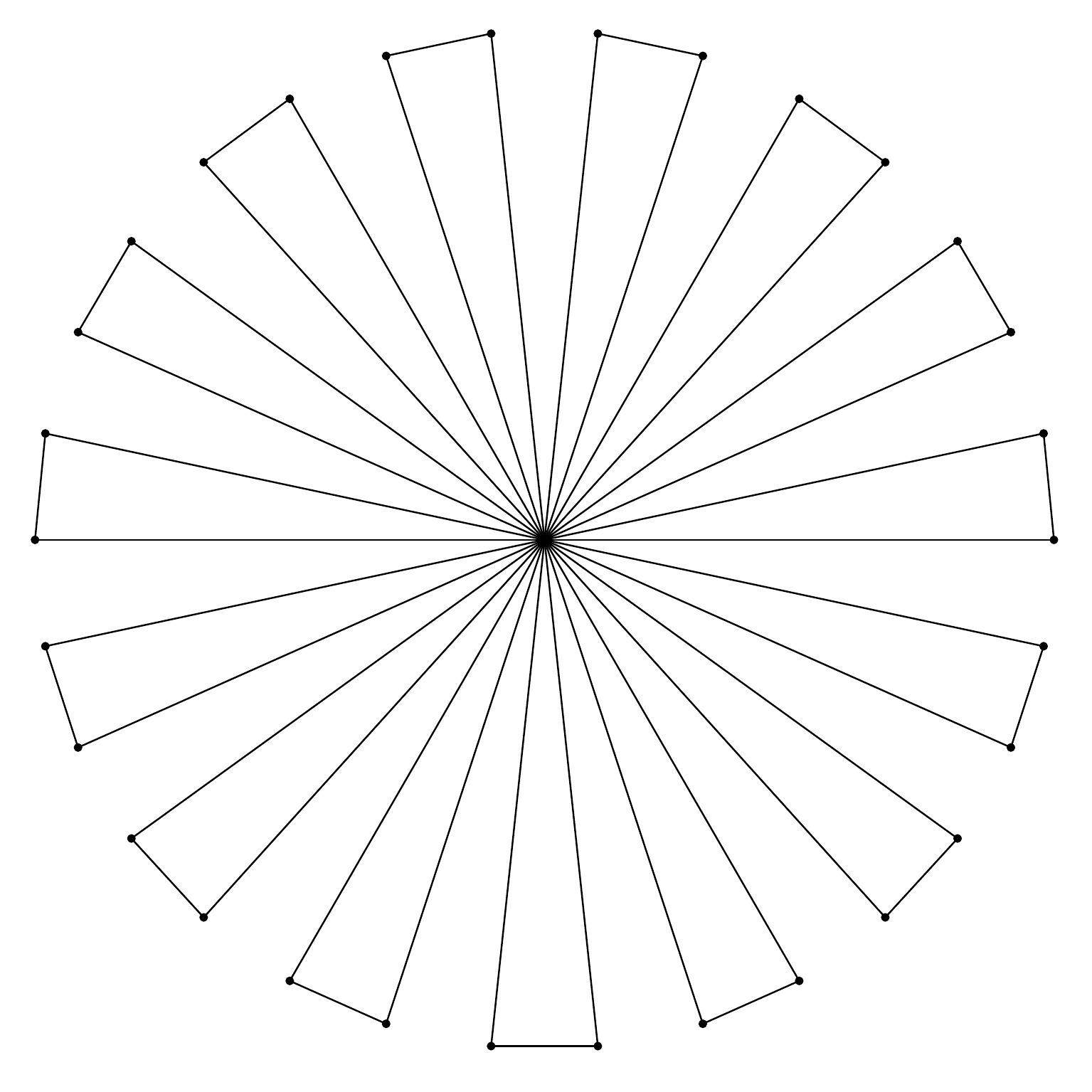} & \includegraphics[width=0.2\textwidth]{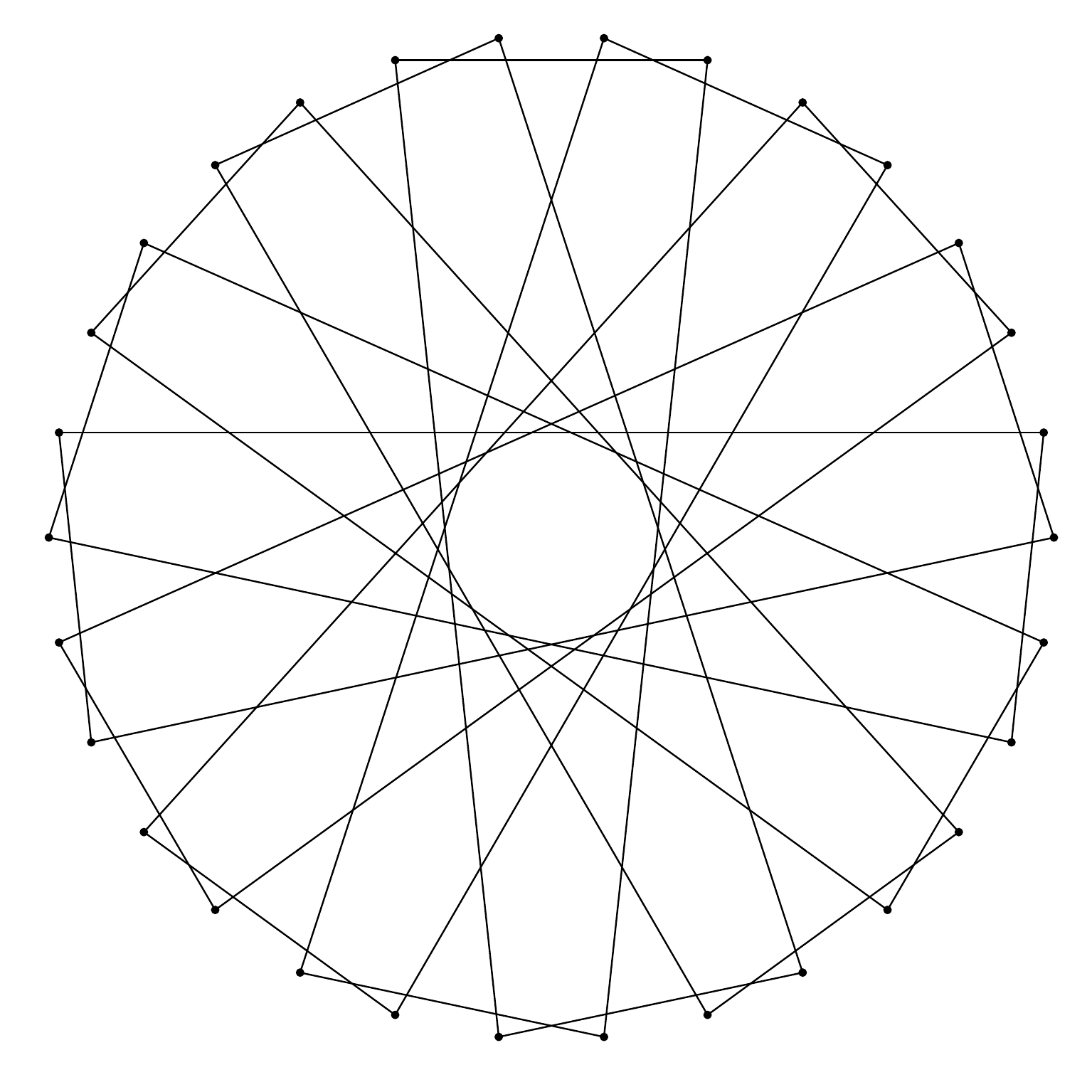}\\
$u=7;a=3;b=11$ & $u=7;a=5;b=9$ & $u=8;a=1;b=15$ & $u=8;a=3;b=13$\\ \hline
\includegraphics[width=0.2\textwidth]{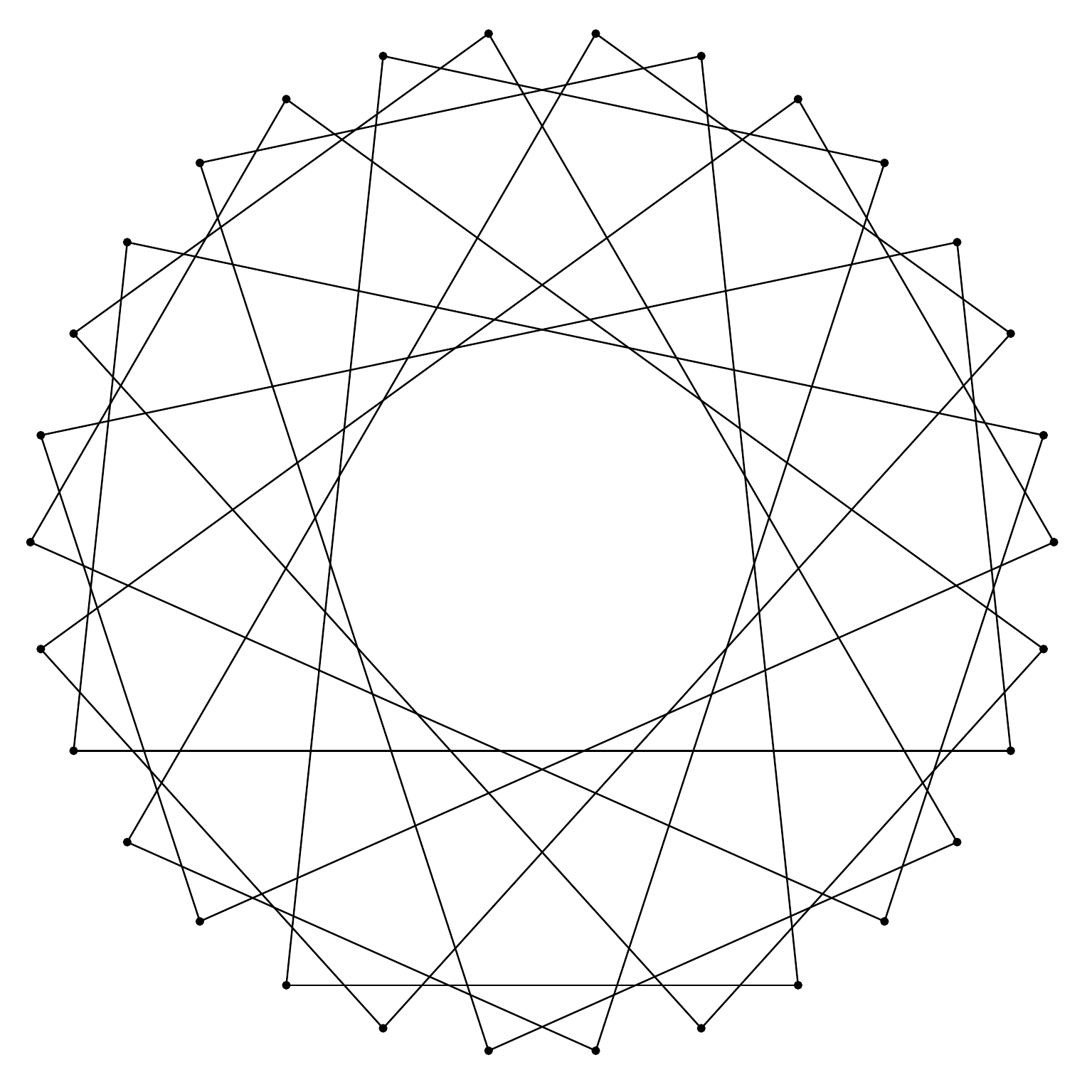} & \includegraphics[width=0.2\textwidth]{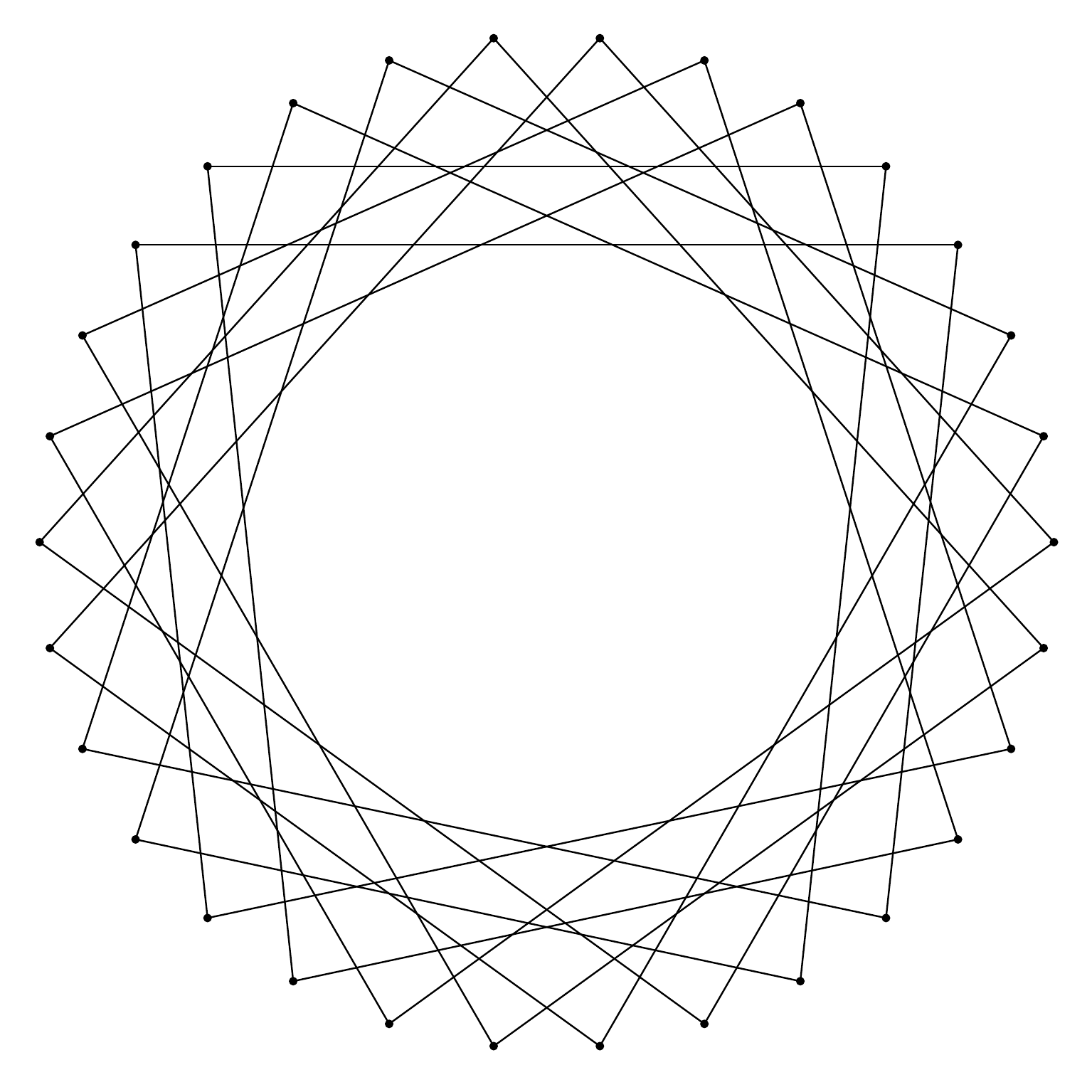} & \includegraphics[width=0.2\textwidth]{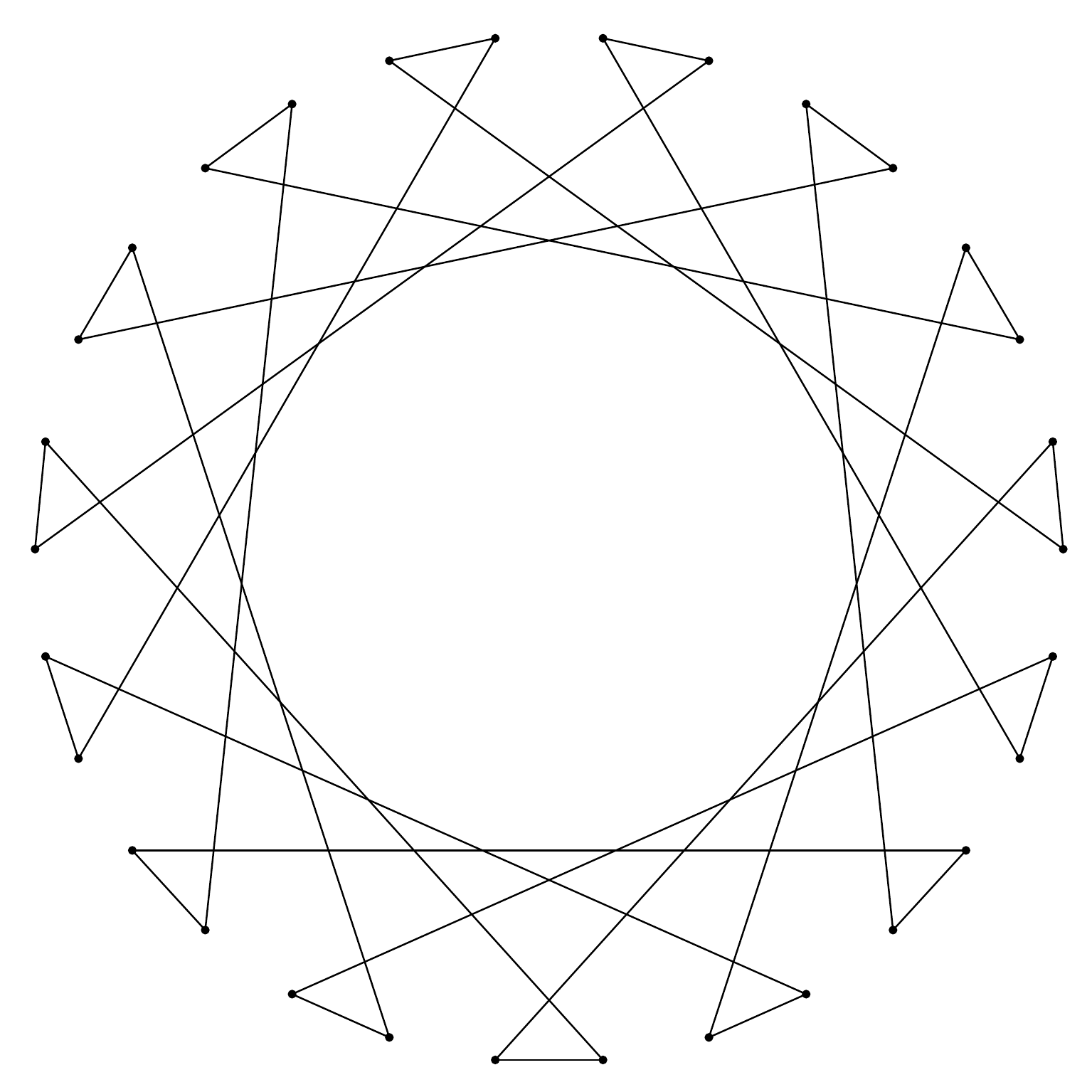} & \includegraphics[width=0.2\textwidth]{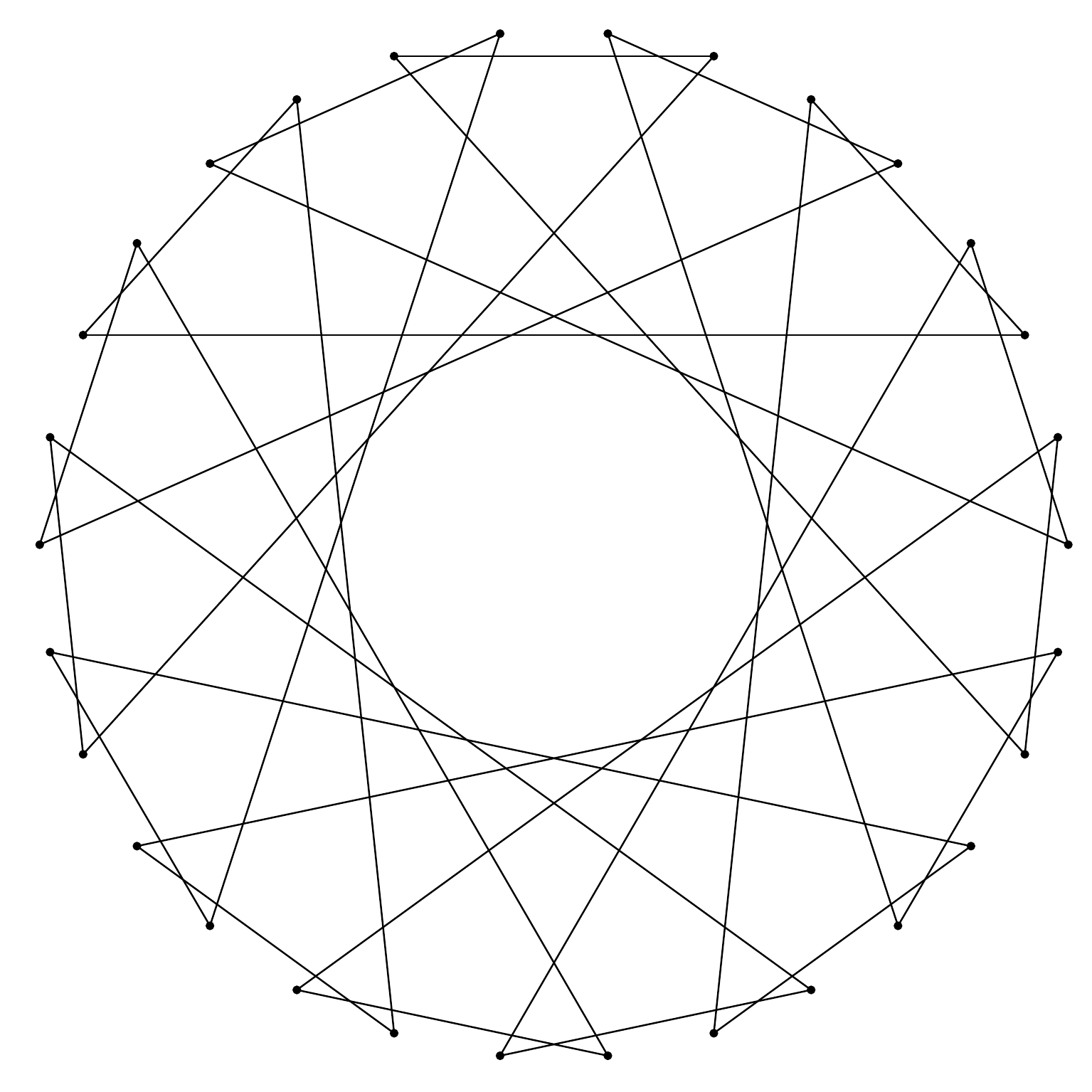}\\
$u=8;a=5;b=11$ & $u=8;a=7;b=9$ & $u=11;a=1;b=21$ & $u=11;a=3;b=19$
\end{tabular}
\caption{The first 12 representatives}
\end{figure}
\newpage
\begin{figure}[!h]
\centering
\begin{tabular}{c | c | c | c}
\includegraphics[width=0.2\textwidth]{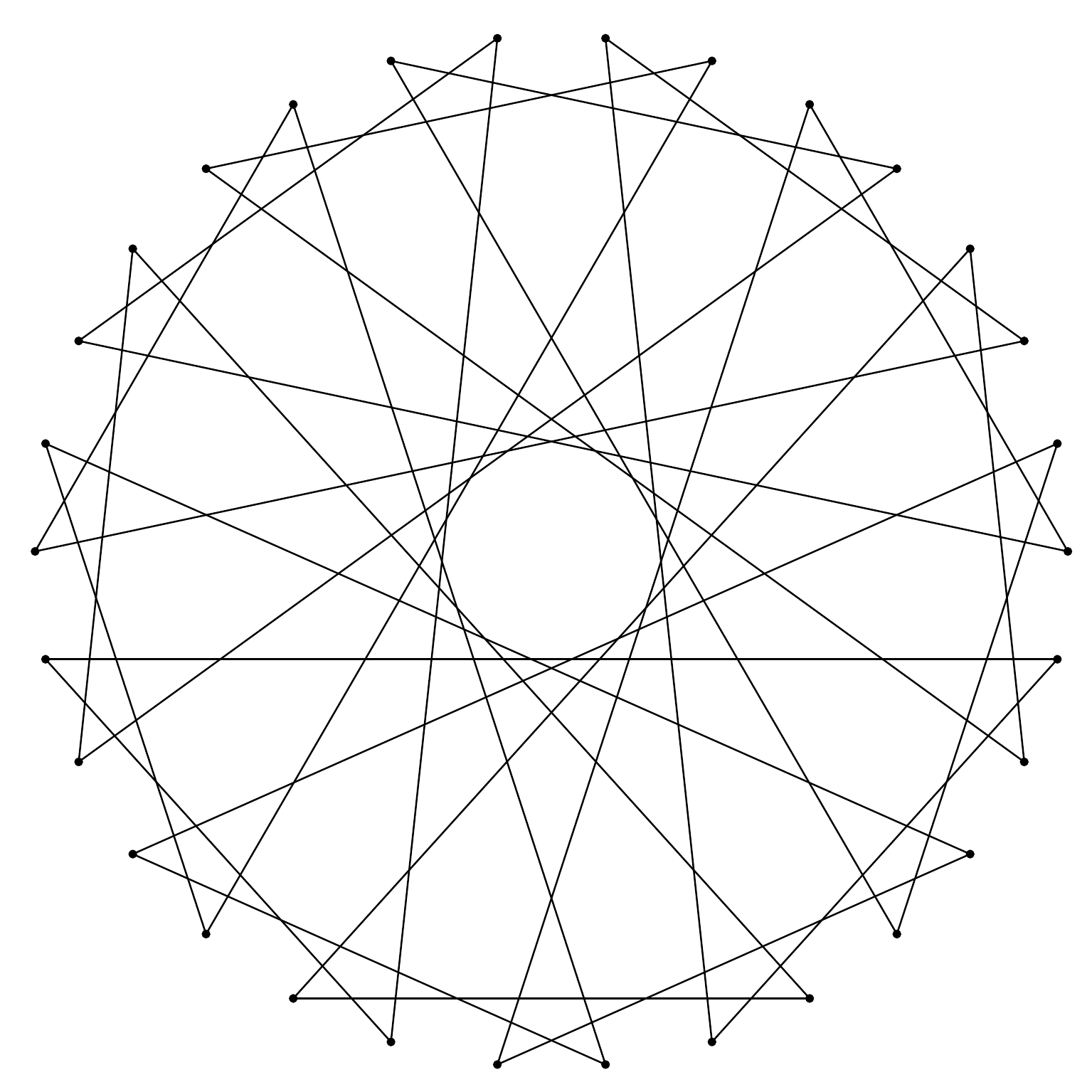} & \includegraphics[width=0.2\textwidth]{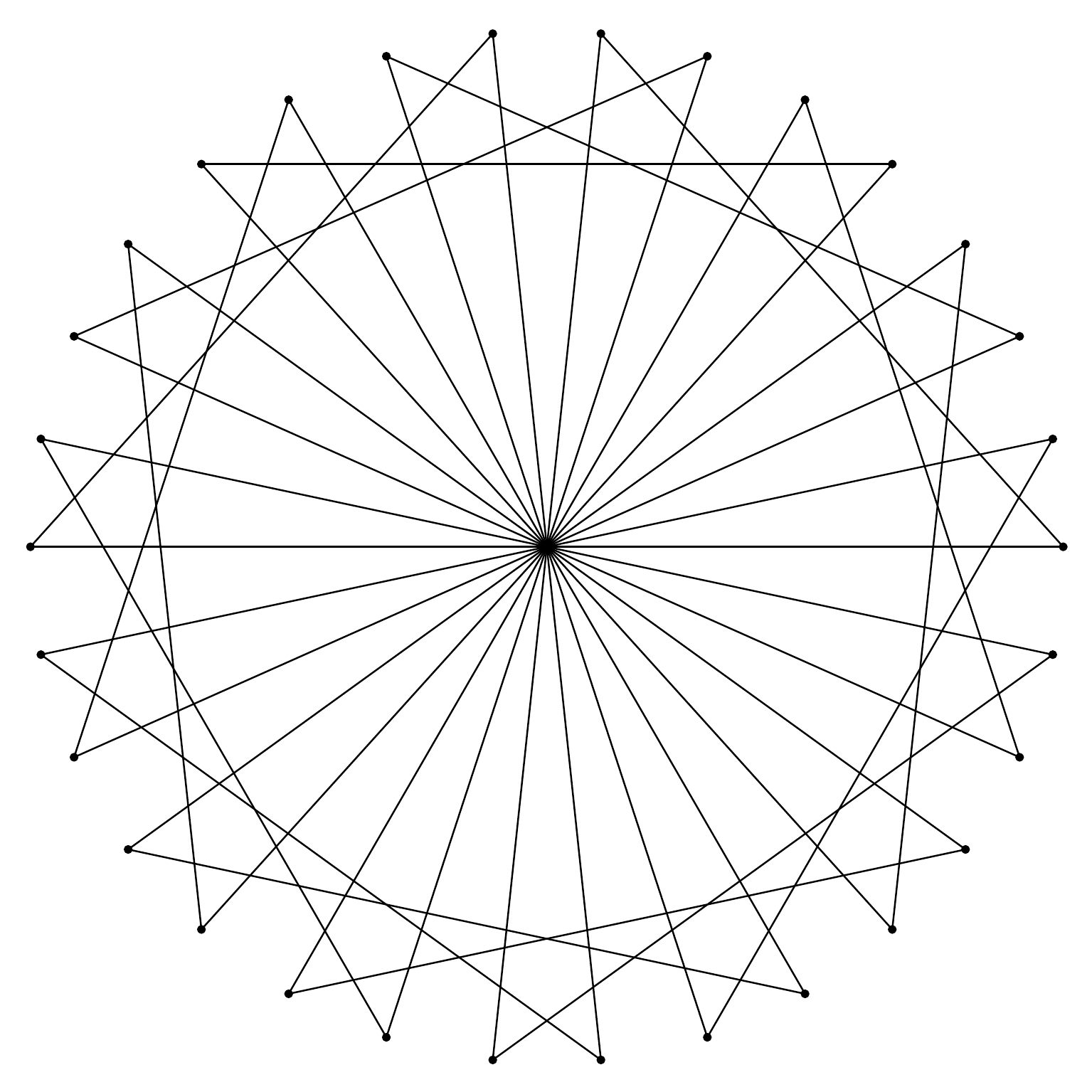} & \includegraphics[width=0.2\textwidth]{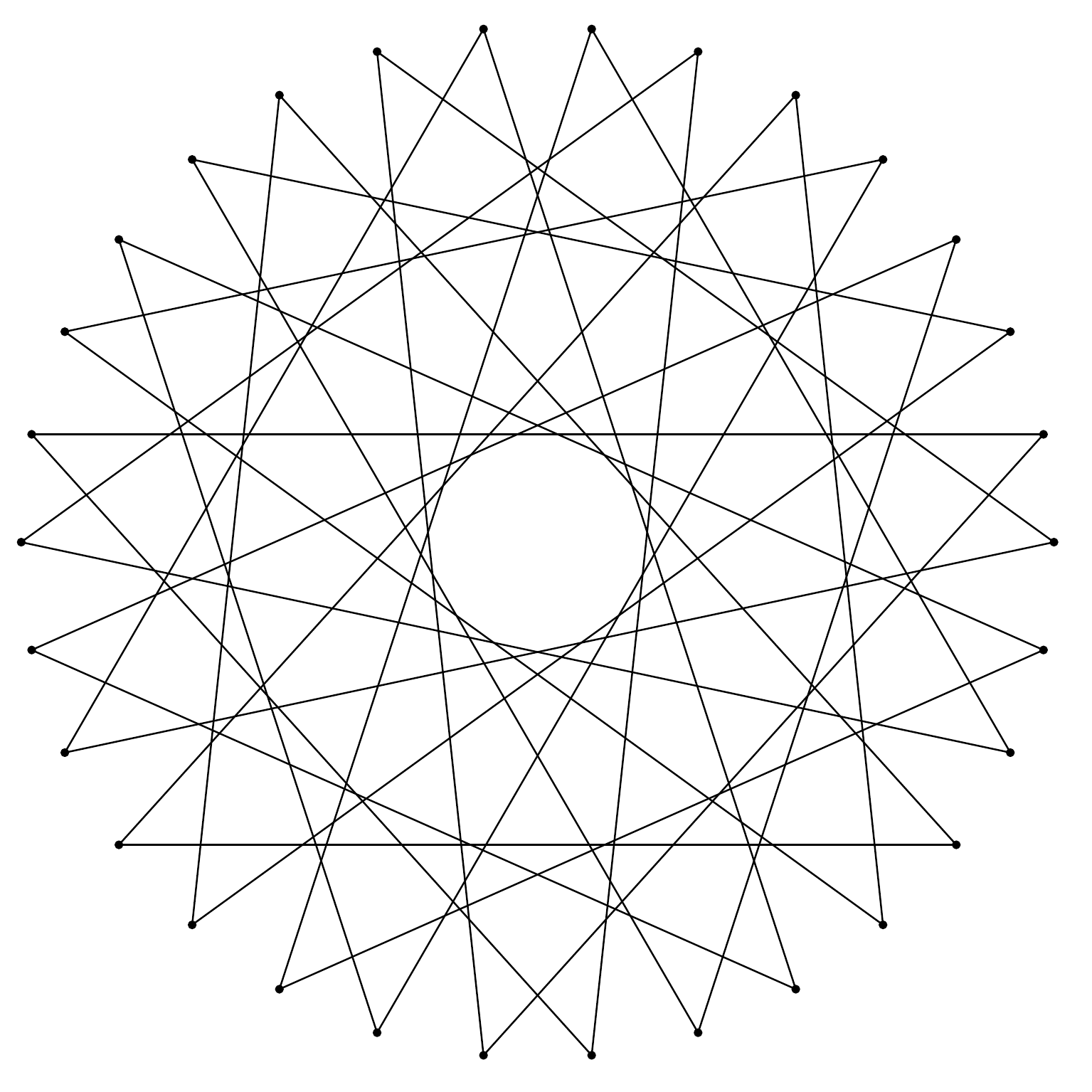} & \includegraphics[width=0.2\textwidth]{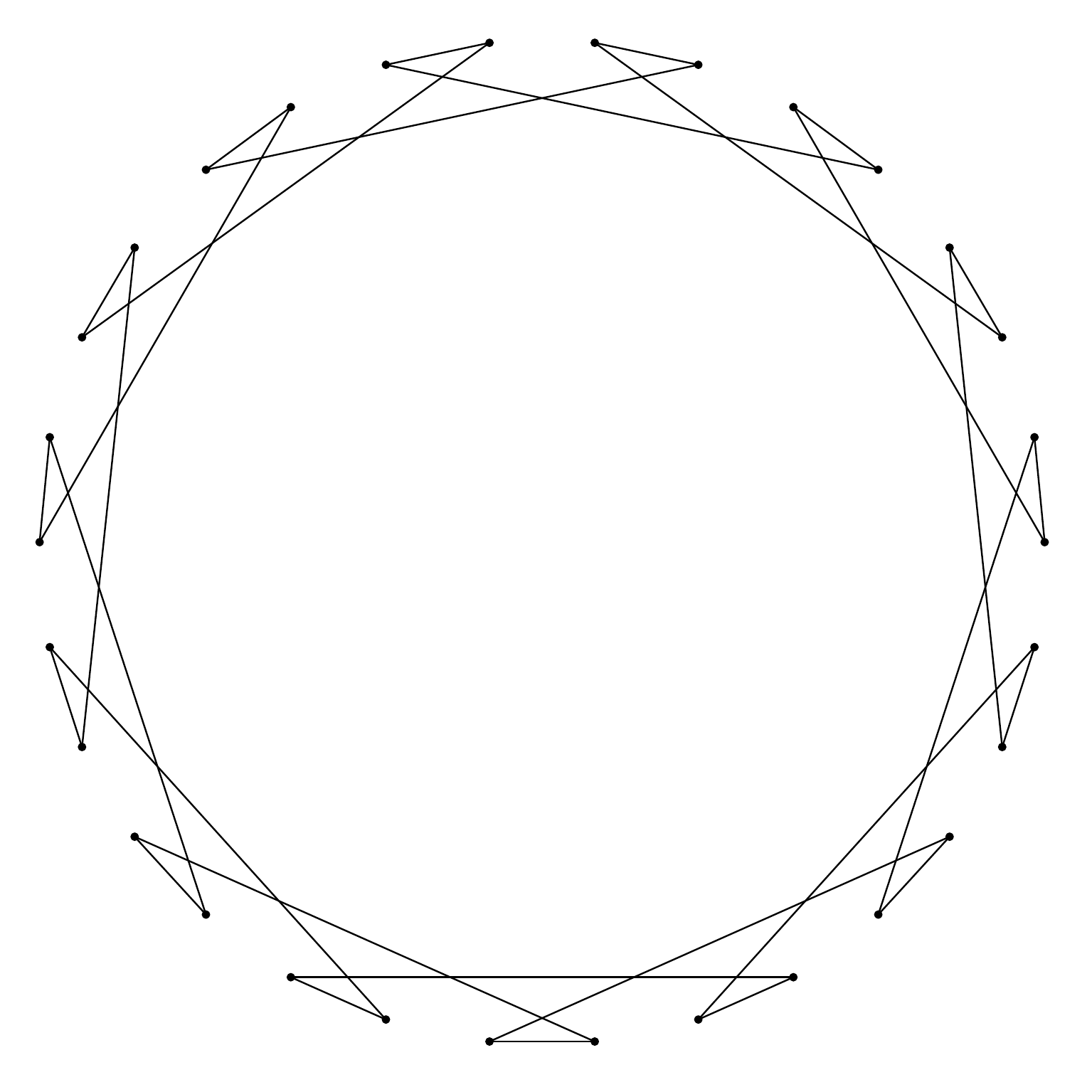}\\
$u=11;a=5;b=17$ & $u=11;a=7;b=15$ & $u=11;a=9;b=13$ & $u=13;a=1;b=25$\\ \hline
\includegraphics[width=0.2\textwidth]{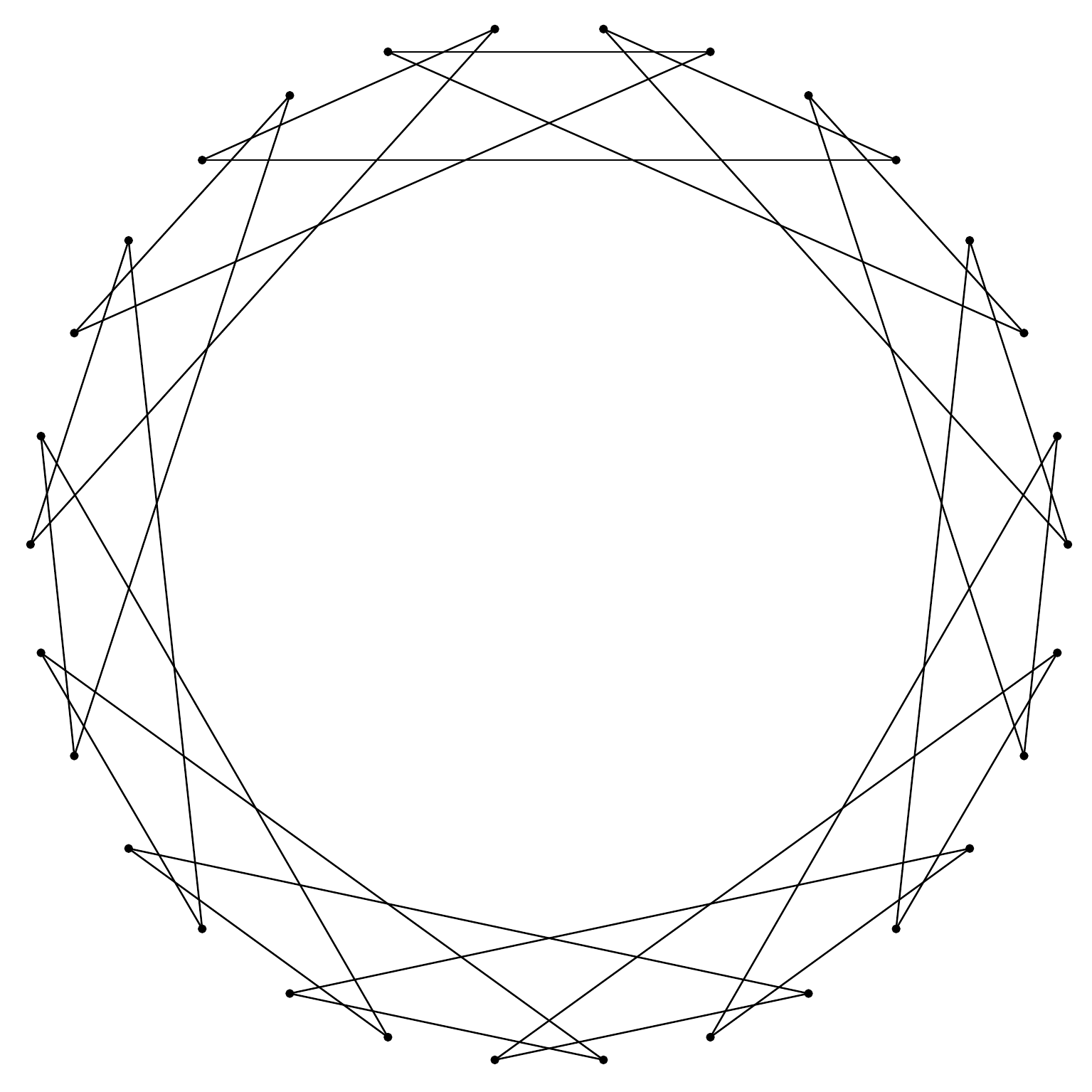} & \includegraphics[width=0.2\textwidth]{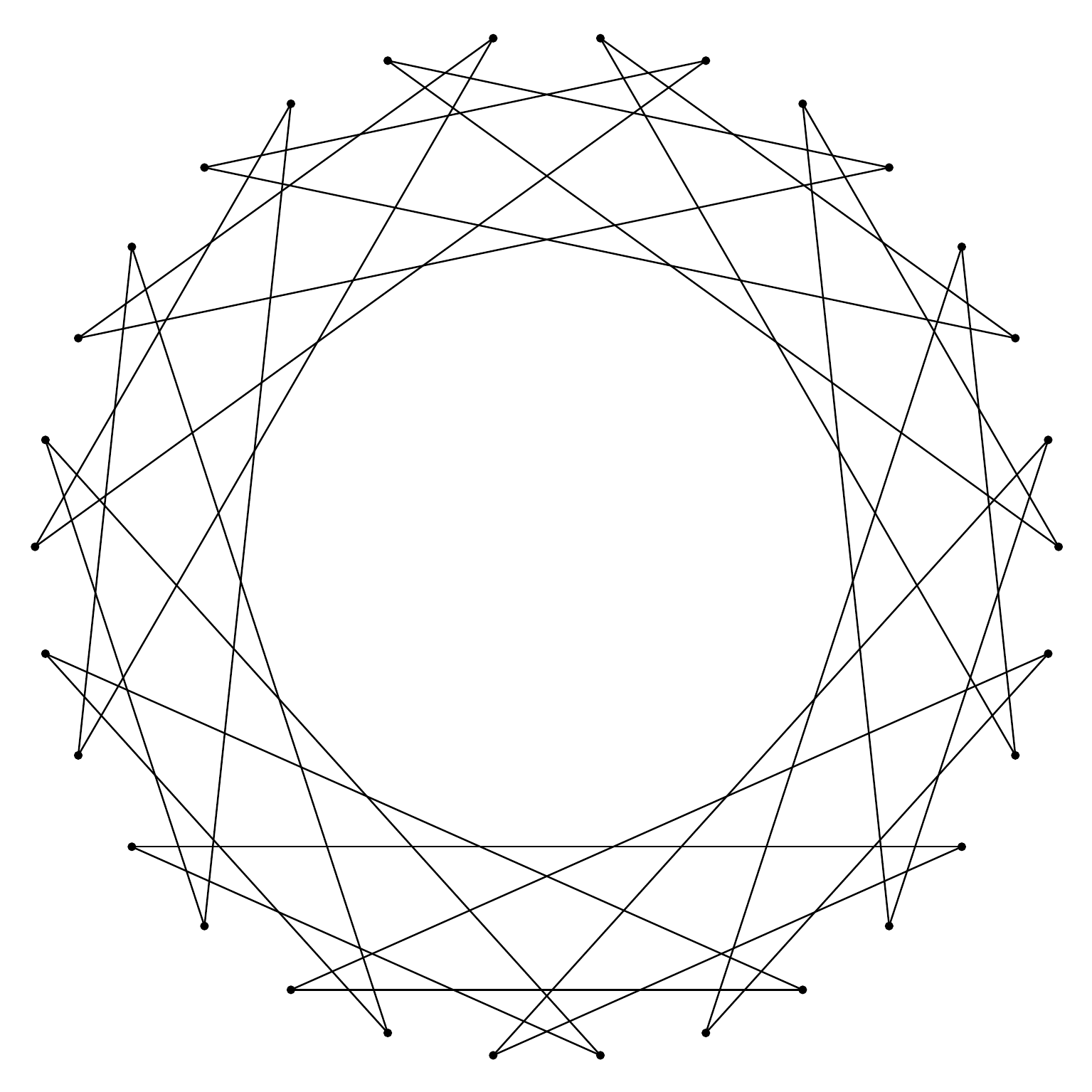} & \includegraphics[width=0.2\textwidth]{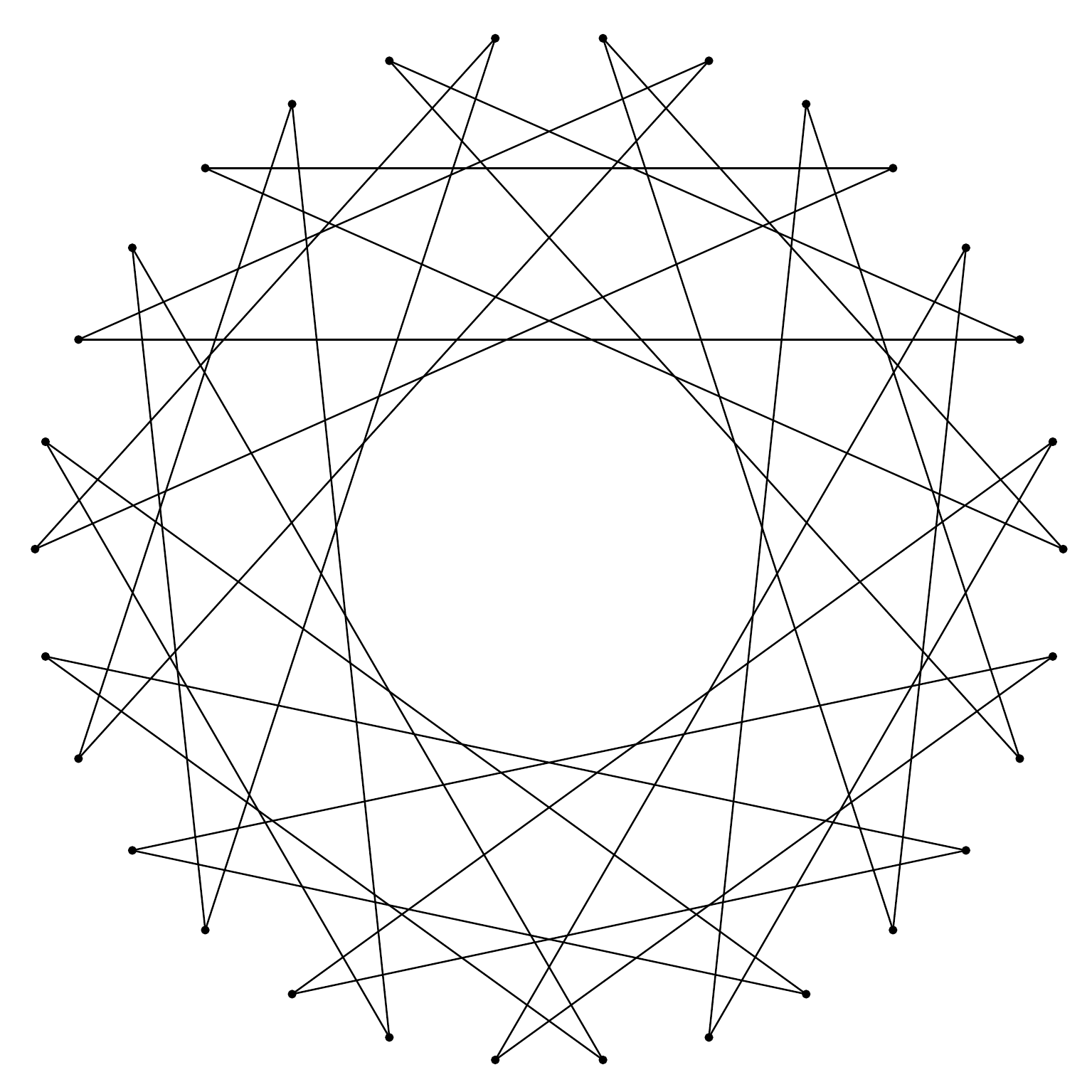} & \includegraphics[width=0.2\textwidth]{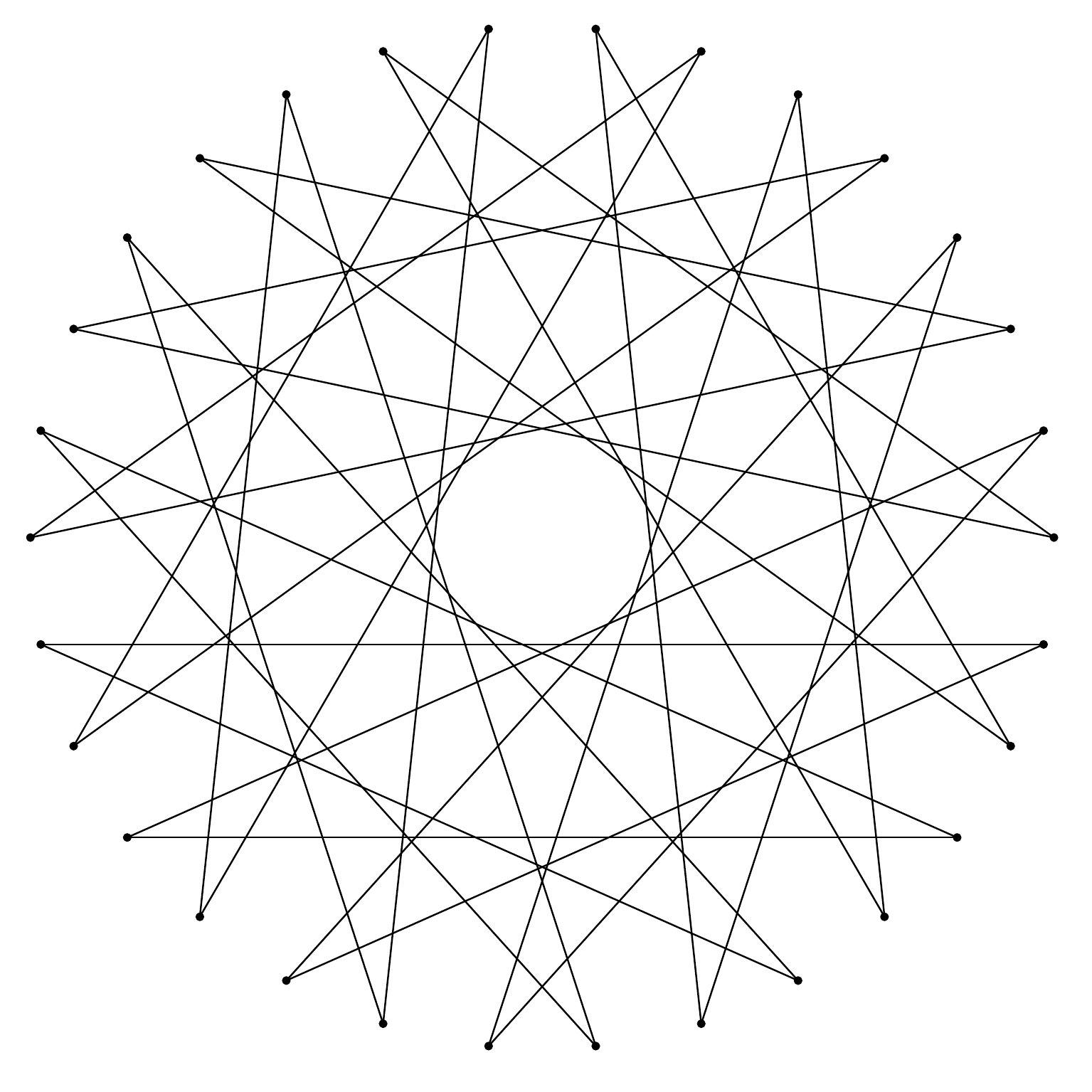}\\
$u=13;a=3;b=23$ & $u=13;a=5;b=21$ & $u=13;a=7;b=19$ & $u=13;a=9;b=17$\\ \hline
\includegraphics[width=0.2\textwidth]{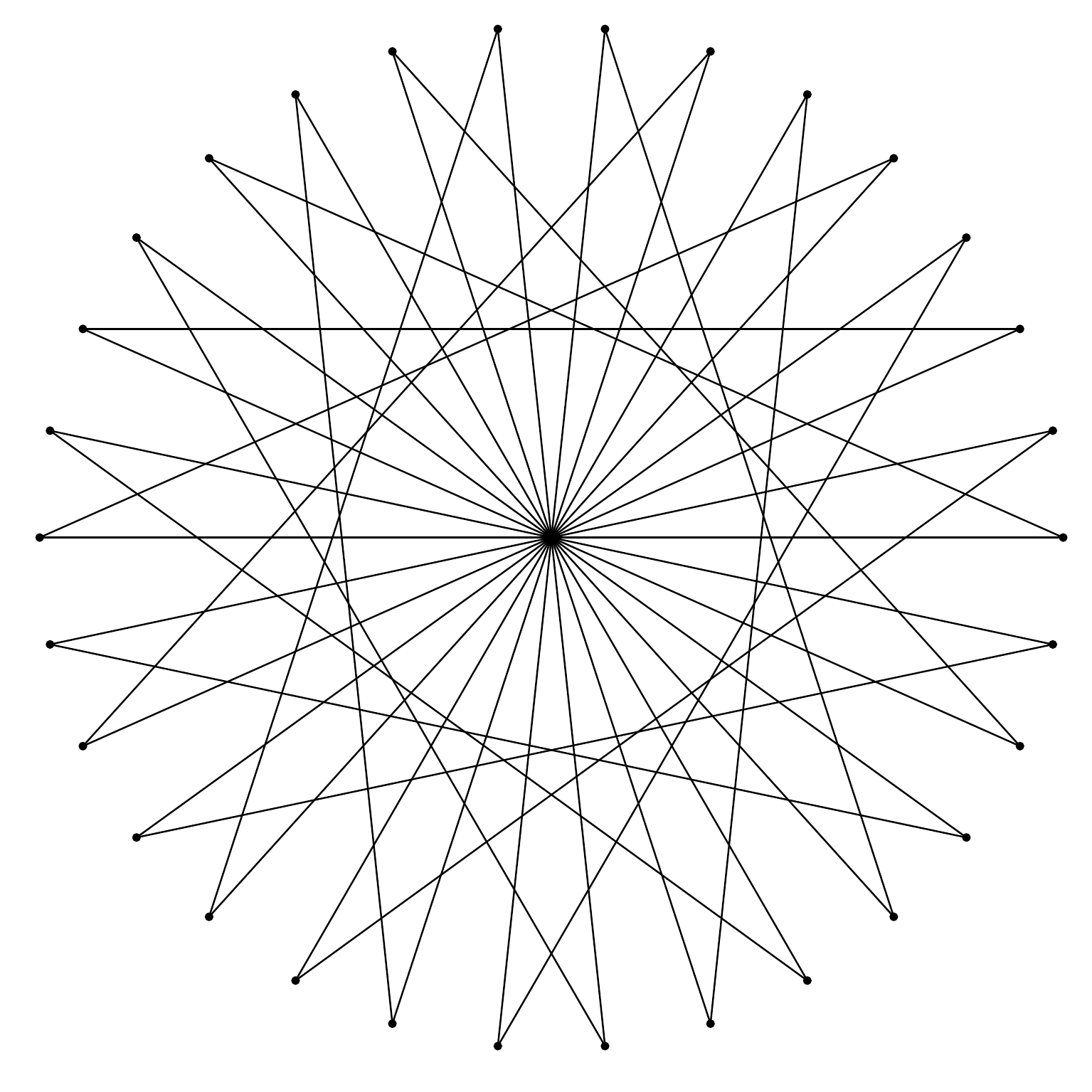} & \includegraphics[width=0.2\textwidth]{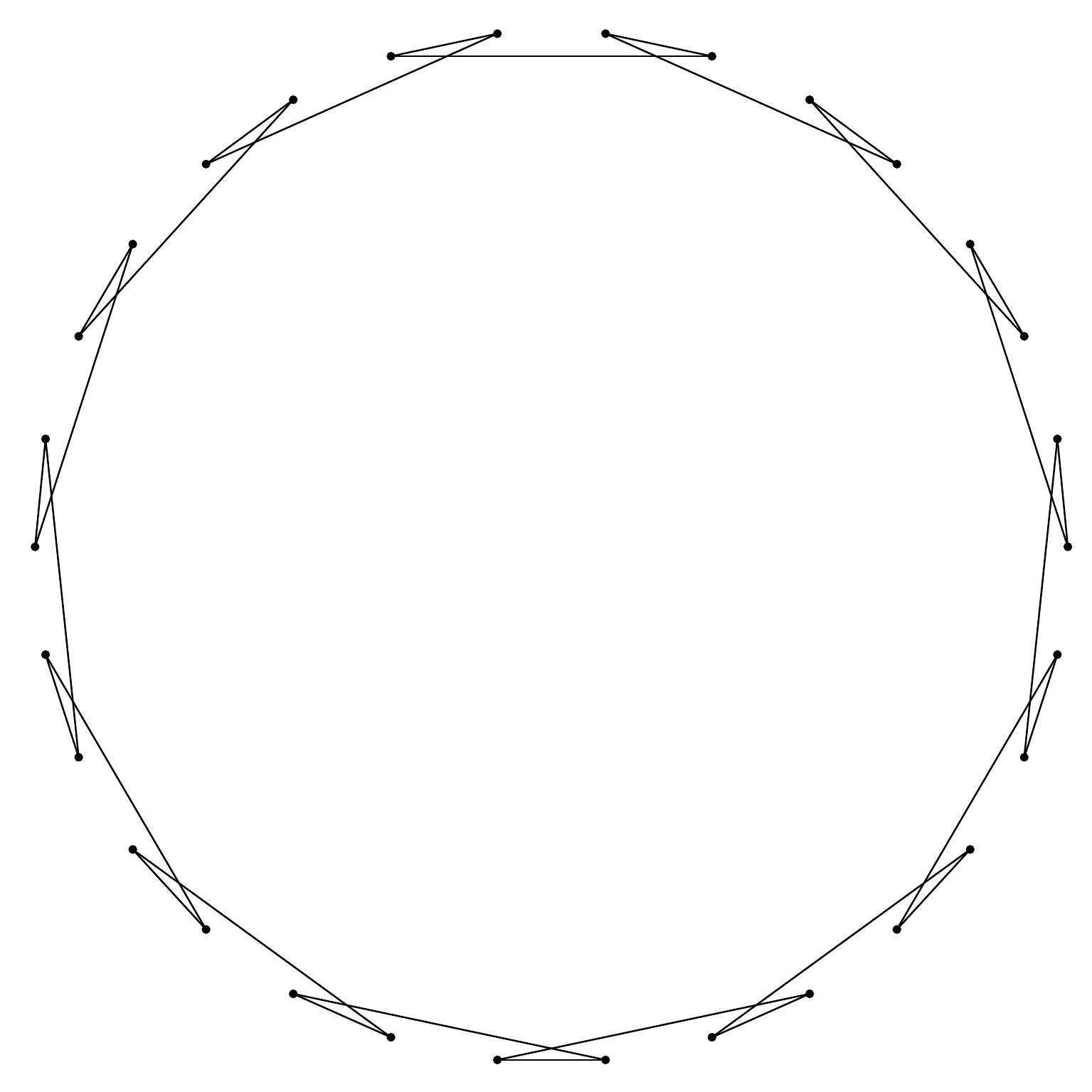} & \includegraphics[width=0.2\textwidth]{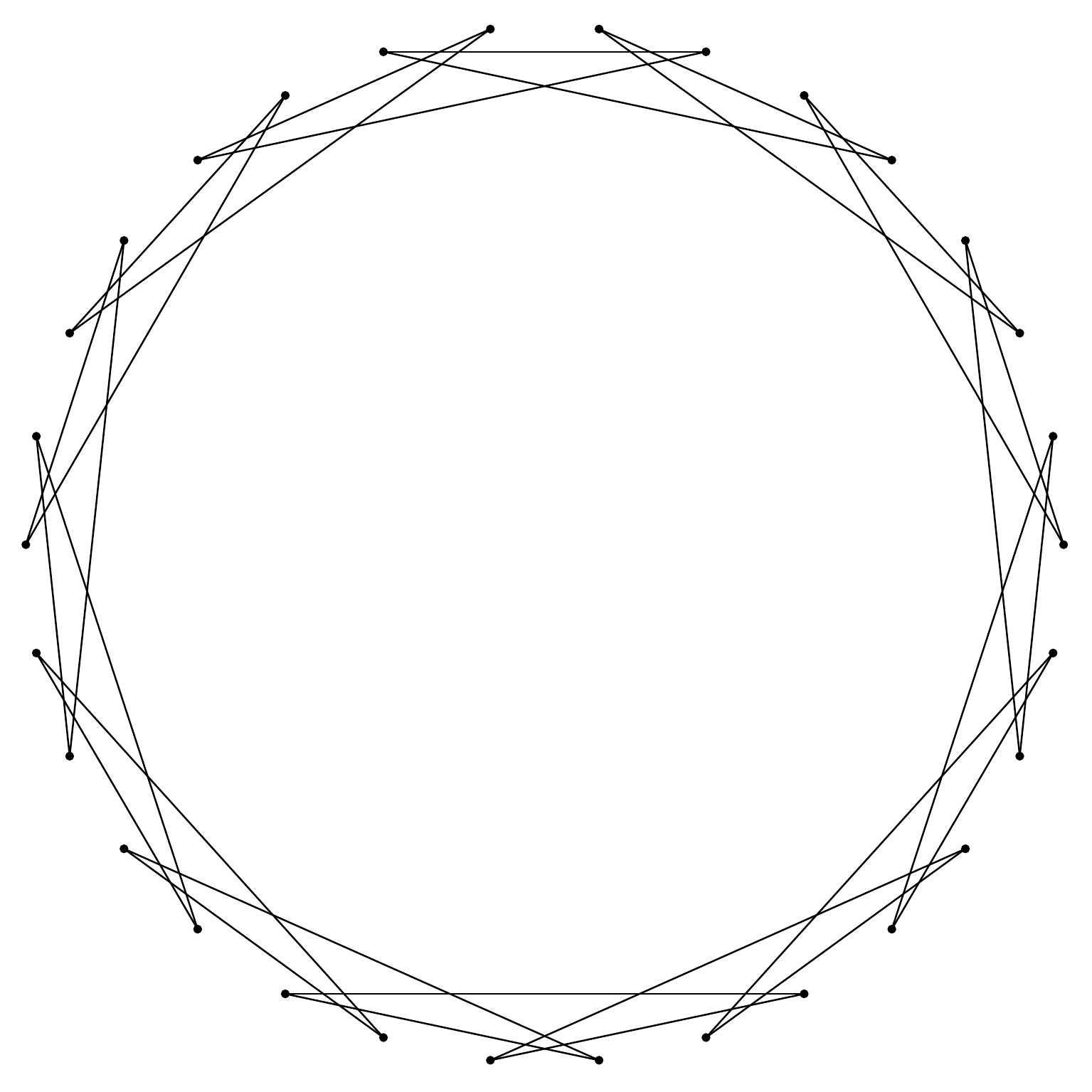} & \includegraphics[width=0.2\textwidth]{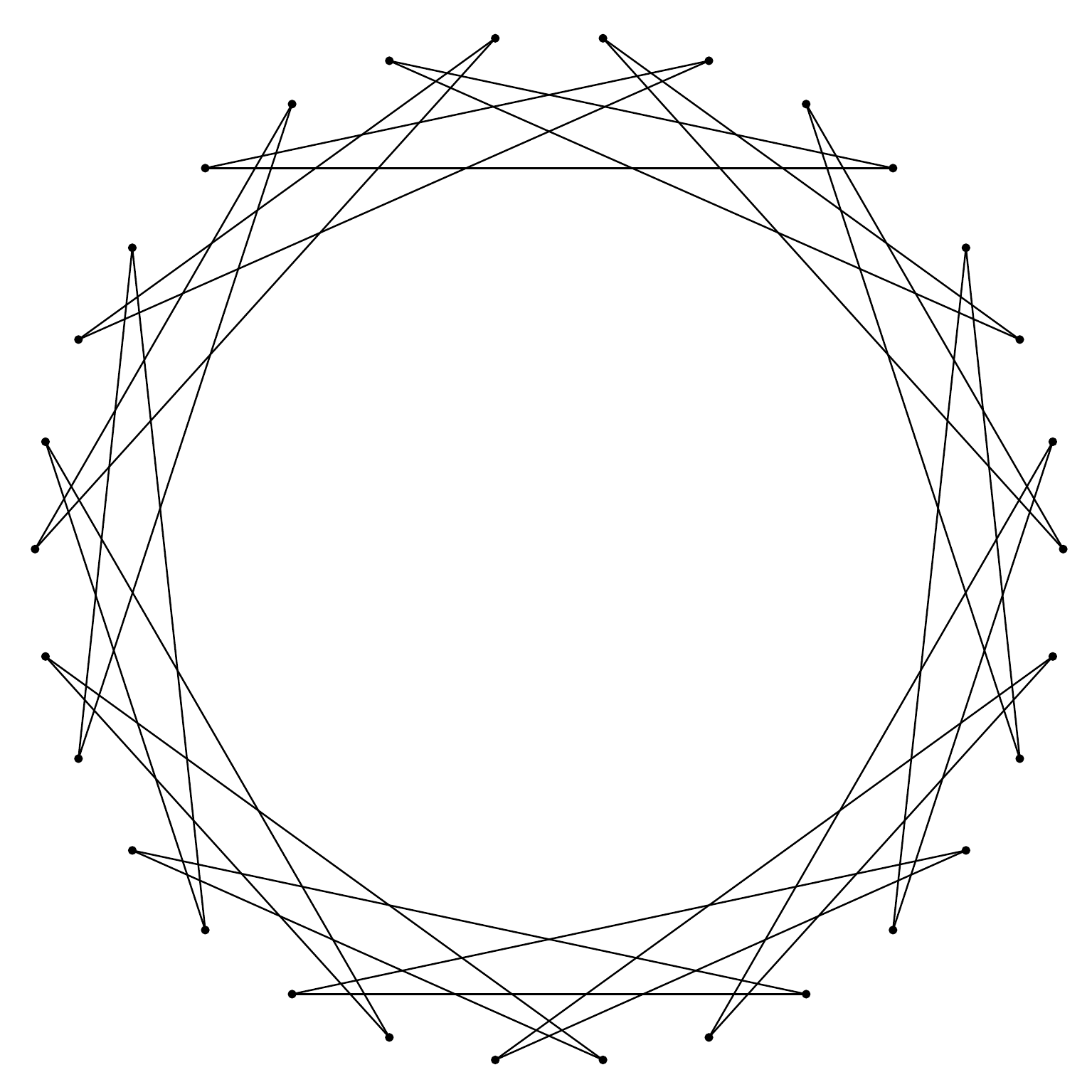}\\
$u=13;a=11;b=15$ & $u=14;a=1;b=27$ & $u=14;a=3;b=25$ & $u=14;a=5;b=23$\\ \hline
\includegraphics[width=0.2\textwidth]{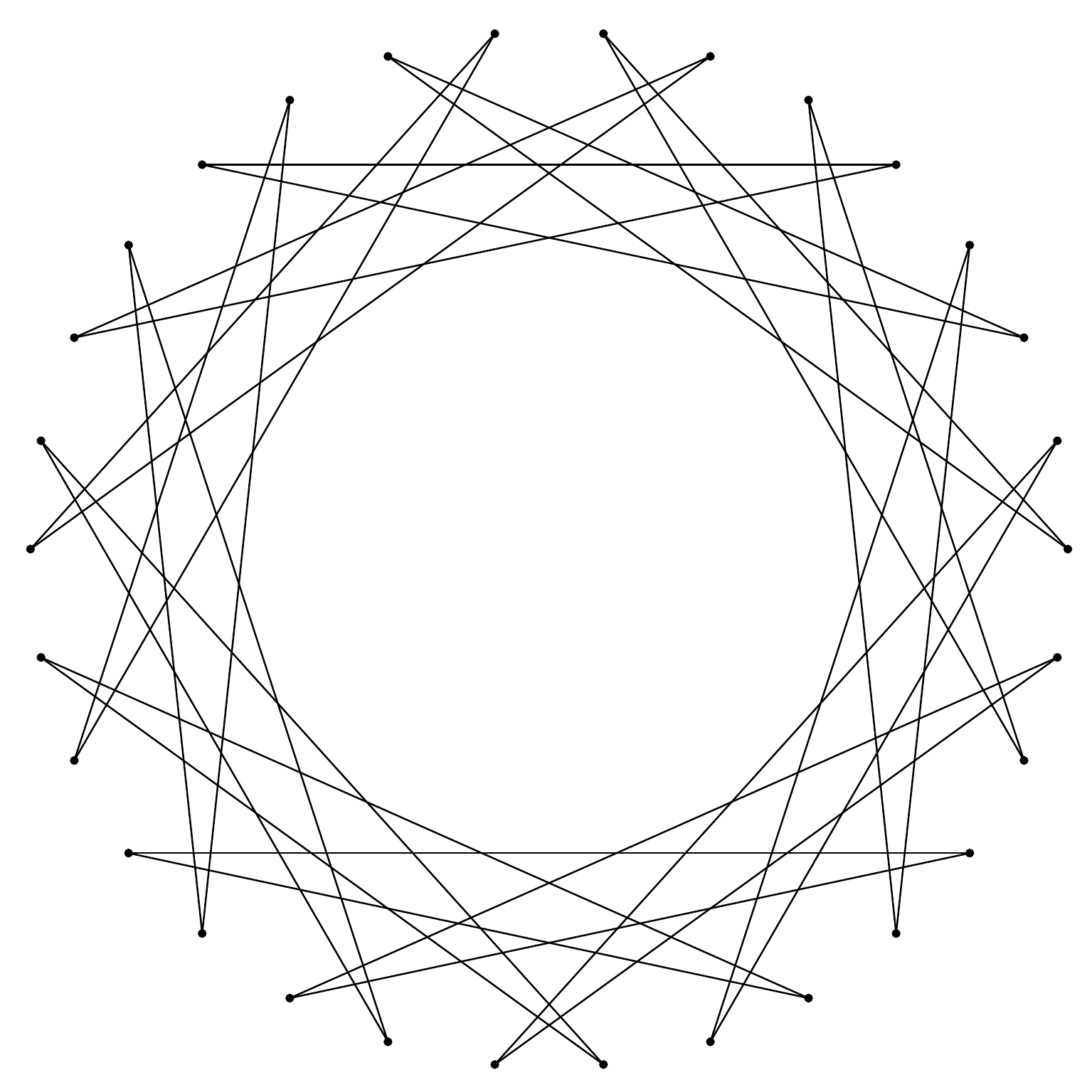} & \includegraphics[width=0.2\textwidth]{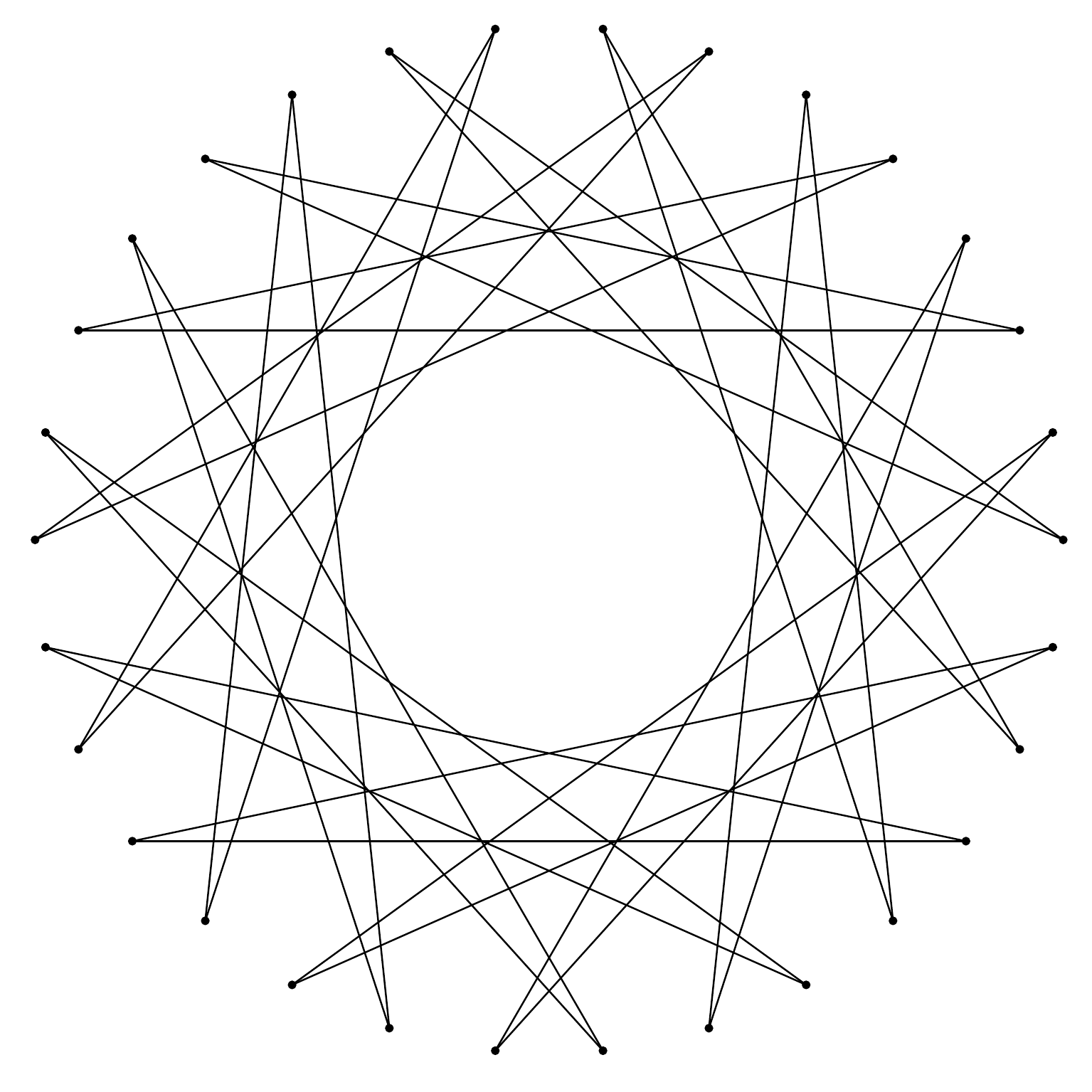} & \includegraphics[width=0.2\textwidth]{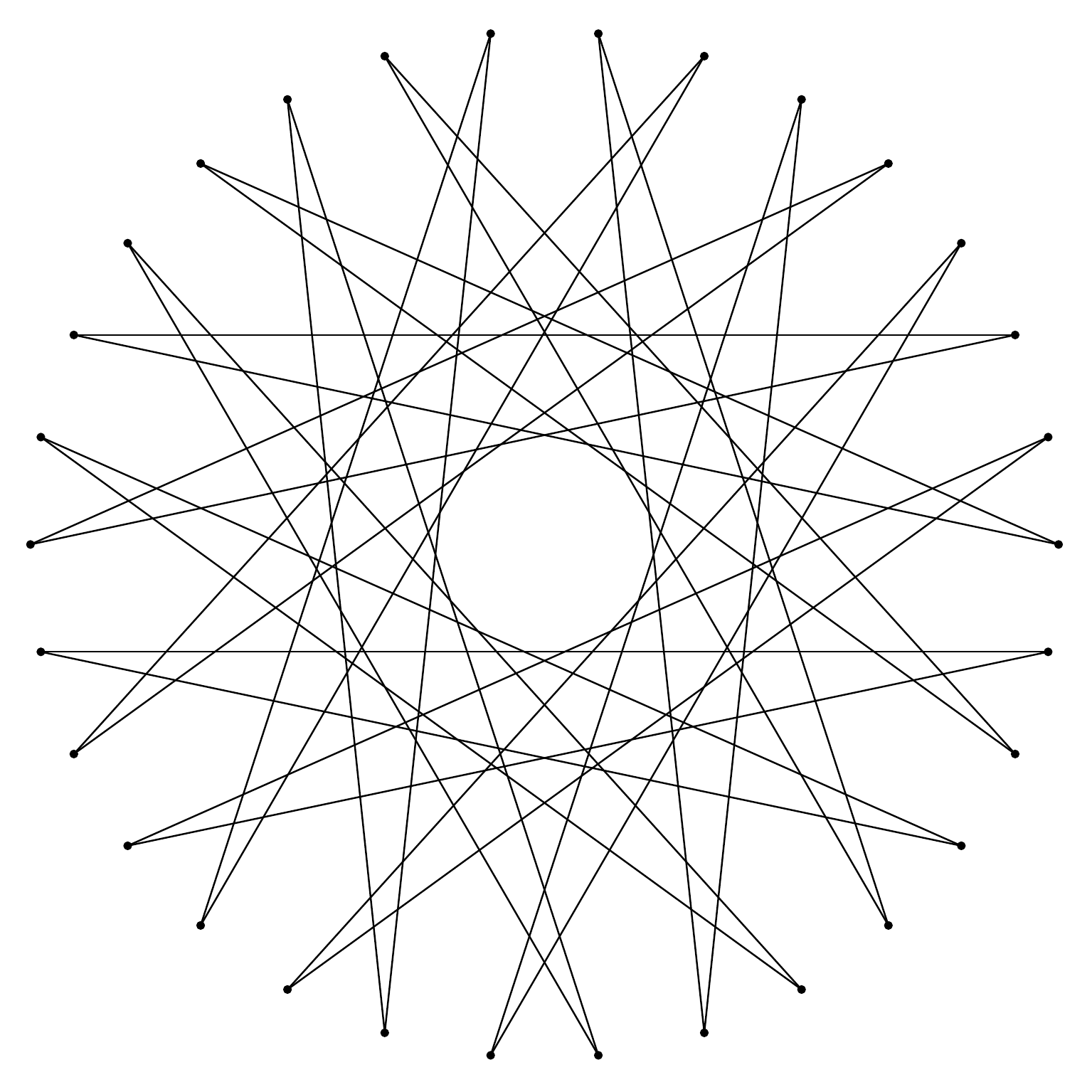} & \includegraphics[width=0.2\textwidth]{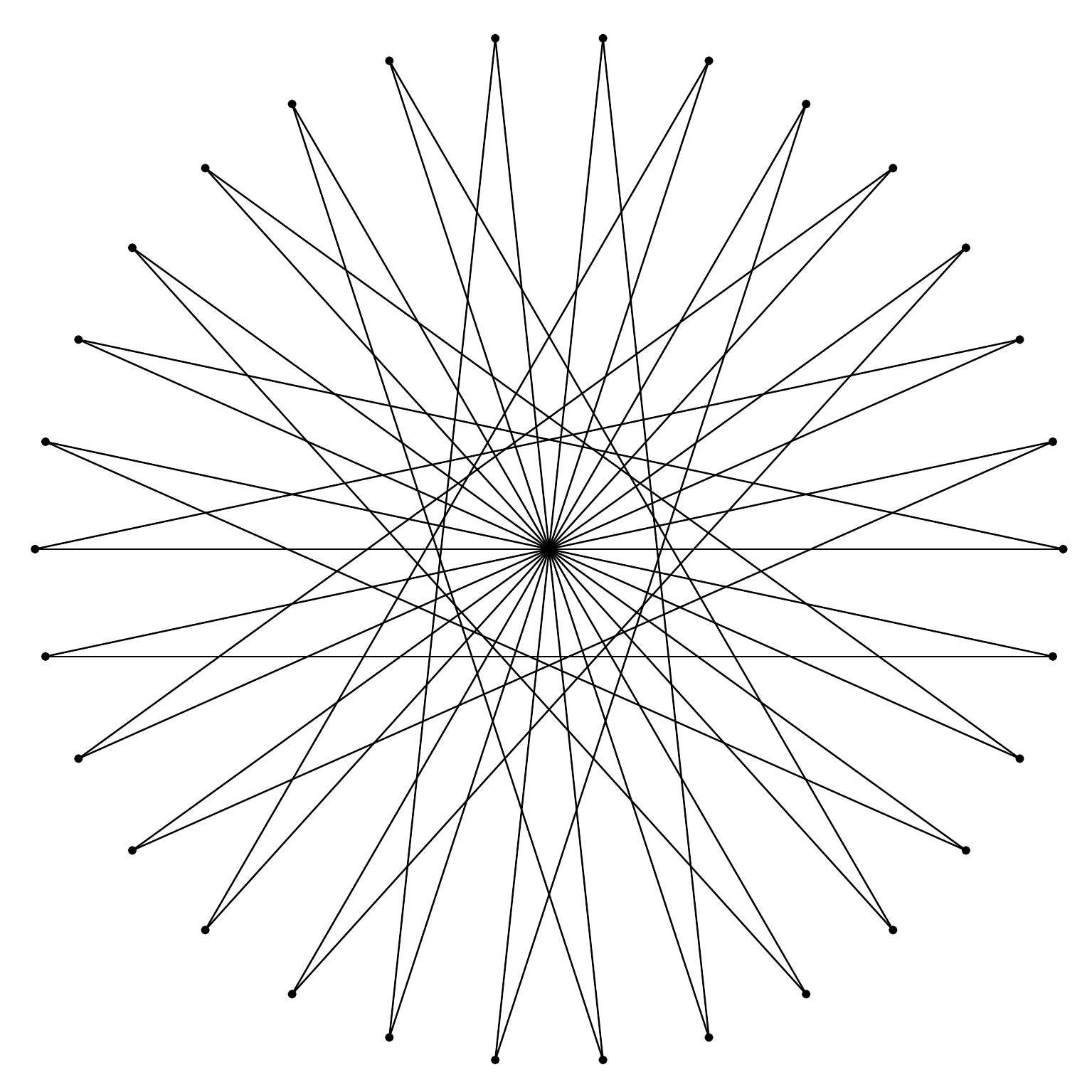}\\
$u=14;a=7;b=21$ & $u=14;a=9;b=19$ & $u=14;a=11;b=17$ & $u=14;a=13;b=15$
\end{tabular}
\caption{The other 16 representatives}
\end{figure}

\subsection{Conclusion: $n=2^{k}$ and the perfect numbers}
\label{subsec:conclusion_2_k_and_the_perfect_numbers}

Let $n=2^{k}$ and $m=2^{k-1}$ with $k\in\mathbb{N}$; $k\geq3$. The possible u-values are the odd numbers from 1 to $n=2^{k-1}-1$. From the main theorem follows immediately:\\

$\vert P_{2^{k-1}}(2^k) \vert=\sum\limits_{u\equiv 1 mod 2=1}^{2^{k-1}-1} \frac{u-1}{2}=\sum\limits_{i=0}^{2^{k-2}-1}i=\dfrac{\left(1+2^{k-2}-1\right)\cdot\left(2^{k-2}-1\right)}{2}=\underline{\underline{2^{k-3}\cdot\left(2^{k-2}-1 \right)}}$.\\
 
This formula contains the perfect numbers. After Euler and Euclid, every perfect number is of the form $ 2^{k-3}\cdot\left(2^{k-2}-1 \right)$ for $k>4$. Thus, we obtain, that the number $\vert P_{2^{k-1}}(2^k) \vert$ is a perfect number, if $2^{k-2}-1$ is a prime number.
\begin{table}[!h]
\begin{center}
\begin{tabular}{| c | c | c | c | c | c |}
\hline
k & n & m & $\vert P_m(n) \vert$ & $2^{k-2}-1$ & perfect\\ \hline
3 & 8 & 4 & 1 & 1 & -\\
4 & 16 & 8 & 6 & 3 & ok\\
5 & 32 & 16 & 28 & 7 & ok\\
6 & 64 & 32 & 120 & 15 & -\\
7 & 128 & 64 & 496 & 31 & ok\\
8 & 256 & 128 & 2016 & 63 & -\\
\hline
\end{tabular}
\caption{Number of equivalence classes of $2^k$-polygons with $2^{k-1}$ axes and perfect numbers}
\label{tab:number_of_equivalence_classes_and_the_perfect_numbers}
\end{center}
\end{table}
\newpage
\subsection{Further conclusions}
\label{subsec:further_conclusions}

\begin{enumerate}
\item Let $p>2$ a prime number. From the main theorem follows in similar way:\\

 $\vert P_p(n) \vert=\vert P_p(2p) \vert=\underline{\underline{\left(\dfrac{n-2}{4} \right)^2}}$. We receive a square number.\\
\item Let $n=2^k\cdot p$ for $k>1$ and $p>2$ prime number. From the main theorem follows in similar way:\\

$\vert P_{2^{k-1}\cdot p}(2^k\cdot p) \vert=\underline{\underline{\dfrac{\left(n-4\right)\cdot\left( n-2^k\right)}{2^5}}}$
\end{enumerate}

\section{Proofs}
\label{sec:proofs}

\subsection{Lemmas}
\label{subsec:lemmas}
Before we begin with the proof of our main theorem, there are a few preparations:\\

\textbf{Lemma 1} Let $a$, $b$ and $n$ be natural even numbers with $a<b\leq n-1$, i.e. $a=2f$, $b=2g$ and $n=2m$ with natural number $f$, $g$ and $m$. If $gcd(a+b,n)=2$ then
\begin{enumerate}\renewcommand{\labelenumi}{\Alph{enumi})}
\item $gcd(f + g, n)=2 \Leftrightarrow  f+g\equiv 0$ mod 2 and $g-f\equiv 0$ mod 2
\item $gcd(f + g, n)=1 \Leftrightarrow f+g\equiv 1$ mod 2 and $g-f\equiv 1$ mod 2
\end{enumerate}

\textbf{Proof of Lemma 1: ($\Rightarrow$)}
\begin{enumerate}\renewcommand{\labelenumi}{\Alph{enumi})}
\item If $gcd(f + g, n)=2$, then $f+g$ must be an even number. And also $g-f$ must be an even number.
\item If $gcd(f + g, n)=1$, then $f+g$ must be an odd number. And also $g-f$ must be an odd number.
\end{enumerate} \hfill $(\Rightarrow)\square$

\textbf{Proof of Lemma 1: ($\Leftarrow$)}
\begin{enumerate}\renewcommand{\labelenumi}{\Alph{enumi})}
\item Let $f+g$ be an even number. Then also $g-f$ must be an even number. Because of the requirements $gcd(a+b,n)=2$ and $a=2f$ and $b=2g$ it follows, that $f<a$, $g<b$ and $gcd(f + g, n) \leq 2$. Because in this case $f+g$ is an even number $gcd(f + g, n) = 2$
\item Let $f+g$ be an odd number. Then also $g-f$ must be an odd number. Because of the requirements $gcd(a+b,n)=2$ and $a=2f$ and $b=2g$ it follows, that $f<a$, $g<b$ and $gcd(f + g, n) \leq 2$. Because in this case $f+g$ is an odd number $gcd(f + g, n) = 1$
\end{enumerate} \hfill $(\Leftarrow)\square$\\

\textbf{Lemma 2} Let $a$ and $b$ be natural odd numbers with $a<b$.\\

$\dfrac{a+b}{2}$ is even $\Leftrightarrow \dfrac{a-b}{2}$ is odd.\\

\textbf{Proof of Lemma 2: ($\Rightarrow$)}

Let $s$ be a natural number, such that $\dfrac{a+b}{2}=2s$. \newline $\Rightarrow a+b=4s \Rightarrow a-b=4s-2b\Rightarrow \dfrac{a-b}{2}=2s-b$. \newline Therefore $\dfrac{a-b}{2}$ must be odd. \hfill $(\Rightarrow)\square$

\textbf{Proof of Lemma 2: ($\Leftarrow$)}

Let $d$ be a natural number, such that $\dfrac{a-b}{2}=2d+1$.\newline  $\Rightarrow a-b=4d+2\Rightarrow a+b=4d+2+2b \Rightarrow \dfrac{a+b}{2}=2d+1+b$. \newline Therefore $\dfrac{a+b}{2}$ must be even. \hfill $(\Leftarrow)\square$\newline\\

\textbf{Lemma 3} Let $a$ and $b$ be natural odd numbers with $a<b$.\\

$\dfrac{a+b}{2}$ is odd $\Leftrightarrow \dfrac{a-b}{2}$ is even.\\

\textbf{Proof of Lemma 3: ($\Rightarrow$)}

Let $s$ be a natural number, such that $\dfrac{a+b}{2}=2s+1$. \newline $\Rightarrow a+b=4s+2 \Rightarrow a-b=4s+2-2b\Rightarrow \dfrac{a-b}{2}=2s+1-b$. \newline Therefore $\dfrac{a-b}{2}$ must be even. \hfill $(\Rightarrow)\square$

\textbf{Proof of Lemma 3: ($\Leftarrow$)}

Let $d$ be a natural number, such that $\dfrac{a-b}{2}=2d$. \newline $\Rightarrow a-b=4d\Rightarrow a+b=4d+2b \Rightarrow \dfrac{a+b}{2}=2d+b$. \newline Therefore $\dfrac{a+b}{2}$ must be odd. \hfill $(\Leftarrow)\square$\\

 \subsection{Proof of the main theorem}
\label{subsec:proof_of_the_main_theorem}
Let $n=2m>3$ be an even integer.\\

If a $n$-polygon with $m$ axes $P_m(n)$ gets rotated twice by the angle $\dfrac{2\pi}{n}$, then it changes over into itself. A first rotation by the angle $\dfrac{2\pi}{n}$ changes $P_m(n)$ with the sequence of sides \newline $(e_1, e_2, \ldots, e_i, \ldots, e_{n-2}, e_{n-1}, e_n)$  into $\overline{P_m(n)}$ with the sequence $(e_n, e_1, e_2, \ldots, e_i, \ldots, e_{n-2}, e_{n-1})$. A second rotation of $P_m(n)$ by the angle $\dfrac{2\pi}{n}$ changes the sequence $(e_1, e_2, \ldots, e_i, \ldots, e_{n-2}, e_{n-1}, e_n)$ into the sequence $(e_{n-1}, e_n, e_1, e_2, \ldots, e_i, \ldots, e_{n-2})$, which represents the same polygon $P_m(n)$. After every second rotation we must get back our $n$-polygon with $m$ axes.\\

With this we have proved the first theorem:\newline\textbf{Theorem 1:} $P_m(n)$ is determined by two sides  $a\in \mathbb{N}$ and $b\in \mathbb{N}$ with $1\leq a< b\leq n-1$ and the alternating sequence of its $n$ sides $(a, b, a, b, \ldots, a, b)$. \\

$a = b$ is excluded, because in this case the polynomial $P (n)$ would have $m=n$ symmetry axes and would be either a common polygon or a regular star-polygon. The precondition $n = 2 m$ would not be fulfilled.\\  

Since the $n$-polygon $P_m(n) = (a, b, a, b, \ldots, a, b)$  and the $n$-polygon $\overline{P_m(n)}= (b, a, b, a, \ldots, b, a)$ belong to the same equivalence-class, we set $a<b$ without affecting the generality, because it is about to enumerate the number of the equivalence-classes. To each equivalence class belong the two $n$-polygons $P_m(n) = (a, b, a, b, \ldots, a, b)$ and $\overline{P_m(n)}= (b, a, b, a, \ldots, b, a)$.\\

\textbf{The sums of the sides and the revolutions of the $n$-polygon $P(n)$:}

So that a $n$-polygon $P(n)$ does not close prematurely, that is, before all $n-1$ other vertices are passed, no sums of $2, 3, \ldots, n-1$ consecutive sides may be divisible by $n$, the sum $s_n$ of all $n$ sides on the other hand, it must be a multiple of $n$. So if $s_n = u\cdot n$, the natural number $u$ is the number of the revolutions in the circle made by the polygon $P(n)$ during its construction. Because $s_n=\dfrac{n}{2}\cdot(a+b)$, which is obvious, we get $u=\dfrac{a+b}{2}$\\

\textbf{The sums of the sides and the revolutions of the $n$-polygon with $m$ axes $P_m(n)$:}

Let $P_m(n) = (a, b, a, b, \ldots, a, b)$ be our $n$-Polygon with $m$ axes. We get the following sums of sides:

\begin{table}[!htp]
\centering
\begin{tabular}{| c | c |}
\hline
$i$ & $ s _i$\\ \hline
1 & $a$ \\
2 & $a+b$\\
3 & $2a+b$\\
4 & $2a+2b=2(a+b)$\\
5 & $3a+2b$\\
6 & $3a+3b=3(a+b)$\\
$\cdots$ & $\cdots$\\
$i$ odd & $\dfrac{i+1}{2}\cdot a+\dfrac{i-1}{2}\cdot b$\\
$i$ even & $\dfrac{i}{2}\cdot a+\dfrac{i}{2}\cdot b=\dfrac{i}{2}\cdot(a+b)$\\
$\cdots$ & $\cdots$\\
$n-1$ & $\dfrac{n}{2}\cdot a+\left(\dfrac{n}{2}-1\right)\cdot b$\\
$n$ & $\dfrac{n}{2}\cdot(a+b)$\\
\hline
\end{tabular}
\caption{Sums of the sides of a $n$-polygon with $m$ axes}
\label{tab:sums_of_the_sides_of_a_n-polygon_with_m_axes}
\end{table}
\newpage
With this we have proven the the basic \textbf{Theorem 2: Existence-theorem:}\\

Let $a$ and $b$ be two natural numbers with $1\leq a< b\leq n-1$. Also, let be the sums of the $n$-tuple $(a, b, a, b, \ldots, a, b)$ determined by the above terms for $s_i$ for all $i$ from 1 to $n$.\\
A equivalence-class of polygons $P_m(n)$ is represented by a $n$-tuple $(a, b, a, b, \ldots, a, b)$ if and only if all sums except the sum $s_n$ are incongruent 0 modulo $n$, but the sum $s_n$ is divisible by $n$.

\textbf{Theorem 3:} If $a$ is even and $b$ odd or $b$ even and $a$ odd, $a$ and $b$ do not induce a $n$-polygon $P_m(n)$.\\

\textbf{Proof of Theorem 3:} From the premise follows immediately that the sum $a + b$ is odd and therefore $\dfrac{a+b}{2}=u$ a fractional number, which is forbidden since $u$ indicates the number of revolutions of the $n$-polygon $P_m(n)$. \hfill $\square$

Thus, $a$ and $b$ must either be both even or both odd.\\

\textbf{Theorem 4:} If $gcd(a+b, n)> 2$, $a$ and $b$ do not induce a $n$-polygon $P_m(n)$.\\

\textbf{Proof of Theorem 4:}
\textit{Idea:}  We prove that under the premise of Lemma 3, among the sums $s_2$ to $s_{n-2}$ with even index $i$, there is at least one sum that is divisible by $n$. This means that the $n$-tuple of sides does not represent a $n$-polygon $P_m(n)$ because the closure happens to early.\\

The minimum value $gcd(a + b, n)$ can take is 3.
Set  $i=\dfrac{2n}{gcd(a+b,n)}$ then $i\leq \frac{2n}{3}<n$. This index $i$ is an even natural number, that is smaller than $n$.
Set further $v=\dfrac{a+b}{gcd(a+b,n)}$. $v$ is a natural number, that is smaller than $u=\frac{a+b}{2}$, because the allowed minimum of the denominator of $v$ is 3.\\

We now show, that the equation $s_i=v\cdot n$ is satisfied by the chosen $i$ and $v$.\newline
$s_i=\dfrac{i}{2}\cdot(a+b)=\dfrac{2n\cdot (a+b)}{2 \cdot gcd(a+b,n)}=\underline{\underline{\dfrac{n\cdot (a+b)}{gcd(a+b,n)}}}$ and \newline
$v \cdot n=\dfrac{a+b}{gcd(a+b,n)}\cdot n=\underline{\underline{\dfrac{n\cdot (a+b)}{gcd(a+b,n)}}}$.\\

Therefore it exists a sum with an index smaller than $n$, which is already divisible by $n$, if $gcd(a + b, n)>2$.\hfill $\square$
\newpage
\textbf{Theorem 5:} If $a$ and $b$ are even numbers and $gcd(a+b,n)=2$, then $a$ and $b$ do not induce a $n$-polygon $P_m(n)$.\\

\textbf{Proof of Theorem 5:} \textit{Idea:} We prove that among the sums $s_3$ to $s_{n-1}$ with odd indices, there is at least one sum, which is divisible by $n$. This means that the $n$-tuple does not represent a polygon $P_m(n)$, because it would close too early.\\

Let $v \leq u=\dfrac{a+b}{2}$ be a natural number and $a = 2f$ and $b = 2g$ with natural numbers $f$ and $g$ with $f<g$.\newline 
The equation $s_i=v\cdot n$ for odd $i$ leads to the following linear diophantine equations (*) and (**) with the variables $i$ and $v$:\\

$\dfrac{i+1}{2}\cdot a+\dfrac{i-1}{2}\cdot b=v\cdot n$ \newline
$\dfrac{i}{2}\cdot(a+b)+\dfrac{1}{2}\cdot(a-b)=v\cdot n$ \newline
$(a+b)\cdot i-2n\cdot v=b-a$ \hfill (*)\newline
$2\cdot(f+g)\cdot i-2n\cdot v=2\cdot(g-f)$ \newline
$(f+g)\cdot i+(-1)n\cdot v=g-f$ \hfill (**)\newline

$gcd(f + g,n)= gcd(f + g, -n)$ is the greatest common divisor of the coefficients of variables $i$ and $v$. Lemma 1 garantees that the constant $g-f$ is divisible by this greatest common divisor. Thus, according to the well-known theorem on linear diophantine equations \cite{Bashmakova1974}, the equation must be solvable with natural numbers $i$ and $v$. It remains to show, that a solution pair $(i, v)$ exists such, that the odd number $i$ obeys the inequality $i \leq n-1$ while $v\leq u$. \\

For this we solve the equation (*) to $i$ and estimate $i$ upwards assuming that $v \leq u$.

$(a+b)\cdot i-2n\cdot v=b-a$ \newline
$(a+b)\cdot i=2n\cdot v+b-a$ \newline
$i=\dfrac{2nv}{a+b}+\dfrac{b-a}{a+b}$ \newline
$i=\dfrac{v}{u}\cdot n+\dfrac{b-a}{a+b}$ \\

Because $\dfrac{v}{u}\leq 1$ and $\dfrac{b-a}{a+b}<1$ we get $i<n+1$ and finally $i\leq n-1$.  \hfill $\square$\\

\textbf{Theorem 6:}
If $a$ and $b$ are odd numbers and $gcd(a+b,n)=2$, then $a$ and $b$ always induce a $n$-polygon $P_m(n)$.\\

\textbf{Proof of Theorem 6:} \textit{Idea:} We prove for all odd indices, that the corresponding sums are not divisible by $n$. We prove further for the even indices, that the corresponding sums are for the first time divisible by $n$ for $i=n$ if $gcd(a+b,n)=2$. If $gcd(a+b,n):=d>2$, then an index $i<n$ is given, for which the corresponding sum is already divisible by $n$. The proof is done separately for the odd and the even indices.\\

\textbf{Odd indices $i$:}
$s_i=\dfrac{i+1}{2}\cdot a+\dfrac{i-1}{2}\cdot b=\dfrac{a+b}{2}\cdot i+\dfrac{a-b}{2}$. \newline
If $\dfrac{a+b}{2}$ is even, $\dfrac{a+b}{2}\cdot i$ is also even, and according to Lemma 2 $\dfrac{a-b}{2}$ is odd and therefore the sum is $s_i$ is odd, and therefore $s_i$ is not divisible by $n$ for every odd $i$.\newline
If $\dfrac{a+b}{2}$ is odd, $\dfrac{a+b}{2}\cdot i$ is also odd, and according to Lemma 3 $\dfrac{a-b}{2}$ is even and therefore the sum is $s_i$ is odd, and therefore $s_i$ is not divisible by $n$ for every odd $i$.\\

\textbf{Even indices $i$:}
So let $i = 2, 4, \ldots, n$ and $s_i=\dfrac{i}{2}\cdot (a+b)$. Furthermore, let $gcd(a+b,n)=:d\geq 2$.
We set the index $i=\dfrac{2n}{d}$. Since $d=gcd(a+b,n)$, $n$ surely is divisible by $d$ and $i=\dfrac{2n}{d}$ is an even number. It follows, $s_i=\dfrac{i}{2}\cdot (a+b)=\dfrac{n\cdot(a+b)}{d}$.\\

We set $v:=\dfrac{a+b}{d}$. Since $d=gcd(a+b,n)$, the sum $a+b$ surely is divisible by $d$, i.e. $v$ is a natural number.\\
 
Now if $d>2$, then $i=\dfrac{2n}{d}<\dfrac{2n}{2}<n$, i.e. $i<n$.  And if $d>2$, then $v=\dfrac{a+b}{d}<\dfrac{a+b}{2}=u$, i.e. $v<u$. And finally, if $d>2$, then $s_i=\dfrac{n\cdot(a+b)}{d}<\dfrac{n\cdot(a+b)}{2}=s_n$, i.e. $s_i<s_n$. This too small sum $s_i$ is divisible by $n$. I.e. $a$ and $b$ do not induce a $n$-polygon $P_m(n)$, if $d>2$.\\

Only if $d=2$ a $n$-polygon is induced by $a$ and $b$. Because if $d=2$, we get: $i=\dfrac{2n}{d}=\dfrac{2n}{2}=n$ and $v=\dfrac{a+b}{d}=\dfrac{a+b}{2}=u$ and $s_i=\dfrac{n\cdot(a+b)}{d}=\dfrac{n\cdot(a+b)}{2}=s_n$, i.e only in this case the conditions of the existence-theorem (theorem 2) are fulfilled.\hfill $\square$\\

We combine now the theorems 1 to 6 into a theorem 7, which is already close to our main-theorem. For theorem 7, no further proof is required.\\

\textbf{Theorem 7:}\\

Let $n=2m>3$ be an even integer.\newline
The $n$-tuples of sides $(a, b, a, b, \ldots, a, b)$ represent different equivalence classes of $n$-polygons $P_m(n)$, if and only if $a$ and $b$ have the following four properties:
\begin{enumerate}
\item $a\in \mathbb{N}$ with $a\equiv $1 mod 2,
\item $b\in \mathbb{N}$ with $b\equiv $1 mod 2,
\item $1 \leq a<b \leq n-1$,
\item $gcd\left(a+b,n \right)=2$.
\end{enumerate}\hfill $\square$\\

Considering, that the condition $gcd(a+b,n)=2$ is equivalent to the condition \newline $gcd\left(2\cdot \dfrac{a+b}{2},2 \cdot m\right)=2$, further to $gcd\left(\dfrac{a+b}{2},m\right)=1$ and finally to $gcd\left(u,m\right)=1$ we get our\\

\begin{center}
\textbf{Main-theorem:}\\
\end{center}

Let $n=2m>3$ be an even integer.\newline
The $n$-tuples of sides $(a, b, a, b, \ldots, a, b)$ represent different equivalence classes of $n$-polygons $P_m(n)$, if and only if $a$ and $b$ have the following four properties:
\begin{enumerate}
\item $a\in \mathbb{N}$ with $a\equiv $1 mod 2,
\item $b\in \mathbb{N}$ with $b\equiv $1 mod 2,
\item $1 \leq a<b \leq n-1$,
\item $gcd\left(u,m\right)=1$.
\end{enumerate}

We conclude the proofs with the proof of the

\begin{center}
\textbf{Formula for $\vert P_m(n) \vert$}\\
\end{center}

\begin{center}
$\vert P_m(n) \vert=\sum\limits_{\substack{u\equiv 0 mod 2,\\gcd \left( u,m \right)=1 }} \frac{u}{2}+\sum\limits_{\substack{u\equiv 1 mod 2,\\gcd \left( u,m \right)=1 }} \frac{u-1}{2}$
\end{center}

\textbf{Proof of the formula for $\vert P_m(n) \vert$:}\\

Let $\dfrac{a+b}{2}:= u$ be prime to $m$. For each allowed $u$-value we determine all pairs $(a, b)$, which satisfy the first three properties in the main theorem. We have to distinguish two cases: either $u$ is even or $u$ is odd. Finally, we determine the number of allowed pairs in both cases.\\

\textbf{$u$ is even:}

The allowed $a$-values start with $a_1=1$, the allowed $b$-values start with $b_1=2u-1$. We obtain a strictly monotone sequence $A$ of $a$-values and a strictly monotone decreasing sequence $B$ of $b$-values until $b-a=2$, i.e.: $A=<1, 3, 5, \ldots, a_j>$ and $B=<2u-1, 2u-3, \ldots, 2u-a_j>$.\\

The equation $b_j=2u-a_j=a_j+2$ leads to $a_j=u-1$ and $b_j=u+1$. Therefore, the number $j$ of terms  of both sequences equals to $j=\dfrac{u}{2}$, i.e. there are $j=\dfrac{u}{2}$ allowed pairs $(a,b)$ starting with $(1,2u-1)$ and ending with $(u-1,u+1)$.\\

\textbf{$u$ is odd:}

The allowed $a$-values start with $a_1=1$, the allowed $b$-values start with $b_1=2u-1$. We obtain a strictly monotone sequence $A$ of $a$-values and a strictly monotone decreasing sequence $B$ of $b$-values until $b-a=4$, i.e.: $A=<1, 3, 5, \ldots, a_j>$ and $B=<2u-1, 2u-3, \ldots, 2u-a_j>$. A further pair is not possible, because $a\neq b$ is a condition.\\

The equation $b_j=2u-a_j=a_j+4$ leads to $a_j=u-2$ and $b_j=u+2$. Therefore, the number $j$ of terms of both sequences equals to $j=\dfrac{u-1}{2}$, i.e. there are $j=\dfrac{u-1}{2}$ allowed pairs $(a,b)$ starting with $(1,2u-1)$ and ending with $(u-2,u+2)$.\hfill $\square$\\ 

\listoffigures
\listoftables
\bibliographystyle{alpha}
\bibliography{Literatur}
\end{document}